# Differentiation of integral Mittag-Leffler and integral Wright functions with respect to parameters

Alexander Apelblat[1] · Juan Luis González-Santander[2]

## Abstract

Derivatives with respect to the parameters of the integral Mittag-Leffler function and the integral Wright function, recently introduced by us, are calculated. These derivatives can be expressed in the form of infinite sums of quotients of the digamma and gamma functions. In some particular cases, these infinite sums are calculated in closed-form with the help of MATHEMATICA. However, parameter differentiation reduction formulas are explicitly derived in order to check some of the results given by MATHEMATICA, as well as to provide many other new results. In addition, we present these infinite sums graphically for particular values of the parameters. Finally, new results for parameter derivatives of the Mittag-Leffler and Wright functions are reported in the Appendices.

**Keywords** Integral Mittag-Leffler function · Integral Wright function · Derivative with respect to parameters

**Mathematics Subject Classification** 33E12 · 33C10 · 33C20

## 1 Introduction

The two-parameter Mittag-Leffler function is an entire function defined by the series [16, Eqn. 10.46.3]:

✉ Juan Luis González-Santander
gonzalezmarjuan@uniovi.es

Alexander Apelblat
apelblat@bgu.ac.il

[1] Department of Chemical Engineering, Ben Gurion University of the Neguev, 84105 Beer Sheva, Israel

[2] Department of Mathematics, Universidad de Oviedo, 33007 Oviedo, Asturias, Spain

$$\mathrm{E}_{\alpha,\beta}(z) = \sum_{k=0}^{\infty} \frac{z^k}{\Gamma(\alpha k + \beta)}, \quad \alpha > 0, \tag{1.1}$$

and, likewise, the Wright function is defined as [16, Eqn. 10.46.1]:

$$\mathrm{W}_{\alpha,\beta}(z) = \sum_{k=0}^{\infty} \frac{z^k}{k!\,\Gamma(\alpha k + \beta)}, \quad \alpha > -1, \tag{1.2}$$

where both series converge for all complex arguments, $z \in \mathbb{C}$, and real values of the second parameter, $\beta \in \mathbb{R}$. However, here we will restrict our study for real arguments, $x \in \mathbb{R}$. Also, we will extend the range of the first parameter of the Mittag-Leffler function to $\alpha \geq 0$, except for the case $\alpha = \beta = 0$.

In our recent paper [2], two new special functions were introduced, namely the integral Mittag-Leffler function:

$$\mathrm{Ei}_{\alpha,\beta}(x) = \int_0^x \frac{\mathrm{E}_{\alpha,\beta}(t) - 1/\Gamma(\beta)}{t} dt = \sum_{k=1}^{\infty} \frac{x^k}{k\,\Gamma(\alpha k + \beta)}, \tag{1.3}$$

and the integral Wright function:

$$\mathrm{Wi}_{\alpha,\beta}(x) = \int_0^x \frac{\mathrm{W}_{\alpha,\beta}(t) - 1/\Gamma(\beta)}{t} dt = \sum_{k=1}^{\infty} \frac{x^k}{k\,k!\,\Gamma(\alpha k + \beta)}. \tag{1.4}$$

It is possible to obtain these integral functions and their Laplace transforms in terms of elementary and special functions in closed-form for many particular cases of the parameters $\alpha$ and $\beta$ with the help of MATHEMATICA program [2].

In a previous investigation [1], the differentiation of the Mittag-Leffler with respect to the parameters was discussed using some integral representations and applying the Laplace approach. In this investigation, we mainly discuss the differentiation of the integral Mittag-Leffler and the integral Wright functions with respect to the parameters. In both investigations, infinite sums of quotients of the digamma and gamma functions arise. These infinite sums have substantial value because they are rarely considered in the literature [4, 7, 14], although recently we found the calculation of some interesting cases in closed-form [8]. Furthermore, there are some interesting applications of these infinite sums of quotients of digamma and gamma functions in fractional calculus literature [9, 17]. For example, they recently appear in the investigation of the Volterra-Prabhakar derivative of distributed order [5, 6, 18, 19]. In particular, they arise in the memory function defining integro-differential operator and in the mean square displacement experimentally measured [11].

In the present investigation, some of the infinite sums found can be expressed in closed-form for particular values of parameters with the aid of MATHEMATICA program. In a complementary way, reduction formulas in closed-form for the parameter derivatives of the integral Mittag-Leffler and integral Wright functions can be explicitly proved, as presented in this paper. However, both approaches overlap in some cases and they provide apparently different results, so that we prove their equivalence. In any case, we have numerically checked all the results presented in this paper.

In order to observe the qualitative behavior of the parametric derivatives of these integral functions, we have illustrated this study with its graphical representation. Although the choice of the limits of the variables and the parameters is arbitrary, in all cases it is the same for comparison reasons.

Finally, we present some new reduction formulas for the parameter derivatives of the Mittag-Leffler function and the Wright function in Appendices 5 and 6.

## 2 Differentiation of the integral Mittag-Leffler function with respect to parameters

Direct differentiation of (1.3) with respect to parameter $\alpha$ gives

$$G(\alpha, \beta; x) = \frac{\partial \text{Ei}_{\alpha,\beta}(x)}{\partial \alpha} = -\sum_{k=1}^{\infty} \frac{\psi(\alpha k + \beta)}{\Gamma(\alpha k + \beta)} x^k, \tag{2.1}$$

where $\psi(z)$ denotes the digamma function [16, Eqn. 5.2.2]:

$$\psi(z) = \frac{\Gamma'(z)}{\Gamma(z)}, \quad z \neq 0, -1, -2, \ldots \tag{2.2}$$

It is worth noting the following property of the digamma function:

$$\psi(z+1) = \frac{1}{z} + \psi(z), \tag{2.3}$$

where $\forall n = 0, 1, \ldots$, we have

$$\psi(n+1) = -\gamma + \sum_{k=1}^{n} \frac{1}{k}, \tag{2.4}$$

being $\gamma$ the Euler-Mascheroni constant. Also, differentiation with respect to parameter $\beta$ leads to

$$H(\alpha, \beta; x) = \frac{\partial \text{Ei}_{\alpha,\beta}(x)}{\partial \beta} = -\sum_{k=1}^{\infty} \frac{\psi(\alpha k + \beta)}{k \Gamma(\alpha k + \beta)} x^k, \tag{2.5}$$

where the following relation is satisfied,

$$\frac{\partial \text{Ei}_{\alpha,\beta}(x)}{\partial \beta} = \int_0^x \frac{\partial \text{Ei}_{\alpha,\beta}(t)}{\partial \alpha} \frac{dt}{t}. \tag{2.6}$$

Note that $G(\alpha, \beta; 0) = H(\alpha, \beta; 0) = 0$, so that we will consider $x \neq 0$ in this Section in order to avoid indeterminate expressions.

As mentioned above, sums of infinite series containing the digamma function appear infrequently in mathematical investigations [7, 8, 14], but a rather relatively large collection of them is given in Brychkov's handbook [4]. However, the sums of the quotients of the digamma and gamma functions that arise in the case of the parameter derivatives of the integral Mittag-Leffler function are available in the literature only for some values of the parameters. Furthermore, they are less obtainable in closed-form than the sums arising in the parameter derivatives of the Mittag-Leffler function. Some of the results presented here (Tables 1, 2, and 3) were derived with the help of MATHEMATICA program. These results are unknown in the mathematical literature in tabular form. In order to illustrate the results given in these tables, we present next some parameter derivatives of the integral Mittag-Leffler function which have been calculated applying MATHEMATICA to (2.1).

$$G(0, \beta; x) = \left.\frac{\partial \text{Ei}_{\alpha,\beta}(x)}{\partial \alpha}\right|_{\alpha=0} = \frac{\psi(\beta)}{\Gamma(\beta)} \frac{x}{x-1}, \quad |x| < 1, \tag{2.7}$$

$$G(1, 0; x) = \left.\frac{\partial \text{Ei}_{\alpha,\beta}(x)}{\partial \alpha}\right|_{\alpha=1,\beta=0} = -x e^x [\ln x - \text{Ei}(-x)], \tag{2.8}$$

$$G(2, 0; x) = \left.\frac{\partial \text{Ei}_{\alpha,\beta}(x)}{\partial \alpha}\right|_{\alpha=2,\beta=0}$$
$$= \sqrt{x}\left[\sinh\sqrt{x}\left(\text{Chi}\sqrt{x} - \ln\sqrt{x}\right) - \cosh\sqrt{x}\,\text{Shi}\sqrt{x}\right], \tag{2.9}$$

$$G(2, 1; x) = \left.\frac{\partial \text{Ei}_{\alpha,\beta}(x)}{\partial \alpha}\right|_{\alpha=2,\beta=1}$$
$$= \cosh\sqrt{x}\left[\text{Chi}\sqrt{x} - \ln\sqrt{x}\right] - \text{Shi}\sqrt{x}\sinh\sqrt{x} - \gamma, \tag{2.10}$$

and

$$G(4, 0; x) = \left.\frac{\partial \text{Ei}_{\alpha,\beta}(x)}{\partial \alpha}\right|_{\alpha=4,\beta=0}$$
$$= \frac{x^{1/4}}{8}\left\{4\sinh\left(x^{1/4}\right)\text{Chi}\left(x^{1/4}\right) + \sin\left(x^{1/4}\right)\left[\ln x - 4\,\text{Ci}\left(x^{1/4}\right)\right]\right.$$
$$\left. - 4\cosh\left(x^{1/4}\right)\text{Shi}\left(x^{1/4}\right) + 4\cos\left(x^{1/4}\right)\text{Si}\left(x^{1/4}\right) - \sinh\left(x^{1/4}\right)\ln x\right\}, \tag{2.11}$$

where the sine and cosine integrals are denoted by

$$\text{Si}(x) = \int_0^x \frac{\sin t}{t} dt, \tag{2.12}$$

$$\text{Ci}(x) = \gamma + \ln x - \int_0^x \frac{1 - \cos t}{t} dt, \tag{2.13}$$

and the hyperbolic sine and cosine integrals, as well as the exponential integrals, are defined by

$$\text{Shi}(x) = \int_0^x \frac{\sinh t}{t} dt, \tag{2.14}$$

$$\text{Chi}(x) = \gamma + \ln x - \int_0^x \frac{1 - \cosh t}{t} dt, \tag{2.15}$$

$$-\text{Ei}(-x) = \text{PV} \int_x^\infty \frac{e^{-t}}{t} dt, \quad x > 0, \tag{2.16}$$

$$\text{E}_1(z) = \int_z^\infty \frac{e^{-t}}{t} dt, \quad z \in \mathbb{C}. \tag{2.17}$$

Note that the result given in (2.7) is a consequence of the geometric series. Indeed, according to (2.1), we have for $|x| < 1$

$$\left.\frac{\partial \text{Ei}_{\alpha,\beta}(x)}{\partial \alpha}\right|_{\alpha=0} = -\sum_{k=1}^\infty \frac{\psi(\beta)}{\Gamma(\beta)} x^k = -\frac{\psi(\beta)}{\Gamma(\beta)} x \sum_{k=0}^\infty x^k = \frac{\psi(\beta)}{\Gamma(\beta)} \frac{x}{x-1}. \tag{2.18}$$

Also, (2.8) is a consequence of the following formula found in [14] for $x > 0$:

$$\sum_{k=0}^\infty \frac{\psi(k+1)}{k!} x^k = e^x [\ln x - \text{Ei}(-x)]. \tag{2.19}$$

We can generalize the results given in (2.8), (2.9), and (2.11), as follows.

**Theorem 2.1** $\forall n = 1, 2, \ldots$ and $0 \leq m < n$, $m \in \mathbb{Z}$, the following reduction formula holds true:

$$\left.\frac{\partial \text{Ei}_{\alpha,\beta}(x)}{\partial \alpha}\right|_{\alpha=n,\beta=-m} = -\frac{x^{(m+1)/n}}{n} \sum_{t=0}^{n-1} \exp\left(\frac{2\pi i t (m+1)}{n}\right) e^\xi [\text{E}_1(\xi) + \ln \xi], \tag{2.20}$$

where

$$\xi = \exp\left(\frac{2\pi i t}{n}\right) x^{1/n}. \tag{2.21}$$

**Proof** According to (2.1) and (2.4), we have

$$\left.\frac{\partial \mathrm{Ei}_{\alpha,\beta}(x)}{\partial \alpha}\right|_{\alpha=n,\beta=-m} = -\sum_{s=1}^{\infty} \frac{\psi(ns-m)}{\Gamma(ns-m)} x^s$$

$$= x^{(m+1)/n} \sum_{s=1}^{\infty} \frac{\left(x^{1/n}\right)^{ns-m-1}}{(ns-m-1)!} \left(\gamma - \sum_{k=1}^{ns-m-1} \frac{1}{k}\right). \tag{2.22}$$

Note that

$$\sum_{s=1}^{\infty} f(ns-m) = \sum_{s=1}^{\infty} \Theta_{n,m}(s) f(s), \tag{2.23}$$

where

$$\Theta_{n,m}(s) := \begin{cases} 1, & \mod(s+m,n) = 0, \\ 0, & \mod(s+m,n) \neq 0, \end{cases}$$

$$= \frac{1}{n} \sum_{t=0}^{n-1} \exp\left(\frac{2\pi i t (s+m)}{n}\right). \tag{2.24}$$

Therefore,

$$\sum_{s=1}^{\infty} f(ns-m) = \frac{1}{n} \sum_{s=1}^{\infty} f(s) \sum_{t=0}^{n-1} \exp\left(\frac{2\pi i t (s+m)}{n}\right)$$

$$= \frac{1}{n} \sum_{t=0}^{n-1} \exp\left(\frac{2\pi i t m}{n}\right) \sum_{s=1}^{\infty} \exp\left(\frac{2\pi i t s}{n}\right) f(s). \tag{2.25}$$

Apply (2) to (2.22), to obtain

$$\left.\frac{\partial \mathrm{Ei}_{\alpha,\beta}(x)}{\partial \alpha}\right|_{\alpha=n,\beta=-m} = \frac{x^{(m+1)/n}}{n} \sum_{t=0}^{n-1} \exp\left(\frac{2\pi i t m}{n}\right) F(\xi), \tag{2.26}$$

where $\xi$ is defined in (2.21) and

$$F(\xi) = \sum_{s=1}^{\infty} \frac{\xi^s}{s!} \left(\gamma - \sum_{k=1}^{s} \frac{1}{k}\right) = \gamma e^{\xi} - \underbrace{\sum_{s=1}^{\infty} \frac{\xi^s}{s!} \sum_{k=1}^{s} \frac{1}{k}}_{S(\xi)}. \tag{2.27}$$

Now, consider the *unit step function*

$$\theta(k) := \begin{cases} 1, & k \geq 0, \\ 0, & k < 0, \end{cases} \quad (2.28)$$

to calculate

$$S(\xi) = \sum_{s=1}^{\infty} \frac{\xi^s}{s!} \sum_{k=1}^{s} \frac{\theta(s-k)}{k} = \sum_{k=1}^{\infty} \frac{1}{k} \sum_{s=1}^{\infty} \frac{\xi^s}{s!} \theta(s-k)$$
$$= \sum_{k=1}^{\infty} \frac{1}{k} \sum_{s=k}^{\infty} \frac{\xi^s}{s!} = \sum_{k=1}^{\infty} \frac{1}{k} \left[e^\xi - e_{k-1}(\xi)\right], \quad (2.29)$$

where $e_k(x) = \sum_{n=0}^{k} \frac{x^n}{n!}$ denotes the *exponential polynomial* [15, Eqn.26:12:2]. Apply now the following formula for the *upper incomplete gamma function* $\forall k = 1, 2, \ldots$ [15, Eqn. 45:4:2]

$$\Gamma(k, x) = (k-1)! \, e_{k-1}(x) \, e^{-x}, \quad (2.30)$$

taking into account that $\Gamma(\nu, x) + \gamma(\nu, x) = \Gamma(\nu)$ [15, Eqn.45:0:1], where $\gamma(\nu, x)$ denotes the *lower incomplete gamma function* [15, Eqn.45:3:1]

$$\gamma(\nu, x) = \int_0^x e^{-t} t^{\nu-1} dt, \quad (2.31)$$

to arrive at

$$S(\xi) = e^\xi \sum_{k=1}^{\infty} \frac{\gamma(k, \xi)}{k!} = e^\xi \int_0^\xi e^{-t} \left(\sum_{k=1}^{\infty} \frac{t^{k-1}}{k!}\right) dt$$
$$= e^\xi \int_0^\xi \frac{1-e^{-t}}{t} dt = e^\xi \, \text{Ein}(\xi), \quad (2.32)$$

where $\text{Ein}(\xi)$ denotes the *entire exponential integral* [16, Eqn. 6.2.3]. Finally, substitute back (2.32) and (2.27) in (2.26), and apply the formula [16, Eqn.6.2.4],

$$\text{Ein}(z) = \text{E}_1(z) + \gamma + \ln z, \quad (2.33)$$

to obtain (2.20), as we wanted to prove. □

**Remark 2.2** The results provided by MATHEMATICA given in (2.8), (2.9), and (2.11) seems to be different from the reduction formula (2.20). However, this is only apparent.

For instance, consider (2.9) and apply the defintions given in (2.14 ) and (2.15) for the hyperbolic sine and cosine integrals to obtain

$$\frac{\partial \text{Ei}_{\alpha,\beta}(x)}{\partial \alpha}\bigg|_{\alpha=2,\beta=0}$$
$$= \sqrt{x}\left[\sinh\sqrt{x}\left(\text{Chi}\sqrt{x} - \ln\sqrt{x}\right) - \cosh\sqrt{x}\,\text{Shi}\sqrt{x}\right]$$
$$= \frac{\sqrt{x}}{2}\left[\left(e^{\sqrt{x}} - e^{-\sqrt{x}}\right)\left(\gamma - \int_0^{\sqrt{x}}\frac{1-\cosh t}{t}dt\right) - \left(e^{\sqrt{x}} + e^{-\sqrt{x}}\right)\int_0^{\sqrt{x}}\frac{\sinh t}{t}dt\right]$$
$$= \frac{\sqrt{x}}{2}\left[e^{\sqrt{x}}\left(\gamma - \int_0^{\sqrt{x}}\frac{1-e^{-t}}{t}dt\right) - e^{-\sqrt{x}}\left(\gamma - \int_0^{-\sqrt{x}}\frac{1-e^{-t}}{t}dt\right)\right]$$
$$= -\frac{\sqrt{x}}{2}\left\{e^{\sqrt{x}}\left[\text{Ein}\left(\sqrt{x}\right) - \gamma\right] - e^{-\sqrt{x}}\left[\text{Ein}\left(-\sqrt{x}\right) - \gamma\right]\right\}$$
$$= -\frac{\sqrt{x}}{2}\left\{e^{\sqrt{x}}\left[\text{E}_1\left(\sqrt{x}\right) + \ln\sqrt{x}\right] - e^{-\sqrt{x}}\left[\text{E}_1\left(-\sqrt{x}\right) + \ln\left(-\sqrt{x}\right)\right]\right\}, \quad (2.34)$$

which is equivalent to (2.20) for $n=2$ and $m=0$.

It is worth noting that applying (2.6) to (2.10 ), and with the aid of MATHEMATICA we arrive at

$$\frac{\partial \text{Ei}_{\alpha,\beta}(x)}{\partial \beta}\bigg|_{\alpha=2,\beta=1} = \frac{x}{4}\,{}_3F_4\left(\begin{array}{c}1,1,1\\2,2,2,\frac{3}{2}\end{array}\bigg|\frac{x}{4}\right)$$
$$+ \left[\ln\sqrt{x} - \text{Chi}\left(\sqrt{x}\right)\right]^2 - \text{Shi}^2\left(\sqrt{x}\right) - \gamma^2, \quad (2.35)$$

where ${}_pF_q\left(\begin{array}{c}a_1,\ldots,a_p\\b_1,\ldots,b_q\end{array}\bigg|x\right)$ denotes the generalized hypergeometric function [16, Sect. 16.2]. The rest of the results given in Table 3 for the derivative of the integral Mittag-Leffler function with respect to $\beta$ for integer values of the parameters $\alpha$ and $\beta$ have been obtained applying MATHEMATICA to (2.5).

Figures 1 and 2 present some graphs of the first derivative of the integral Mittag-Leffler function with respect to the first parameter. $\partial \text{Ei}_{\alpha,\beta}(x)/\partial\alpha$ is plotted as a function of variable $x$ and parameter $\alpha$ at fixed value of parameter $\beta=0$ in Fig. 1. $\partial \text{Ei}_{\alpha,\beta}(x)/\partial\alpha$ is plotted as a function of parameter $\alpha$ at fixed values of $x$ and $\beta=1$ in Fig. 2. Many other graphs of this kind have been plotted for other values of parameter $\beta$. In all of these graphs, we have observed that $\partial \text{Ei}_{\alpha,\beta}(x)/\partial\alpha$ tend to zero for large values of parameter $\alpha$.

Figures 3 and 4 present some graphs of the first derivative of the integral Mittag-Leffler function with respect to the second parameter. $\partial \text{Ei}_{\alpha,\beta}(x)/\partial\beta$ is plotted as a function of variable $x$ and parameter $\alpha$ at fixed value of parameter $\beta=0$ in Fig. 3. $\partial \text{Ei}_{\alpha,\beta}(x)/\partial\beta$ is plotted as a function of parameter $\alpha$ at fixed values of $x$ and $\beta=1$ in Fig. 4. Many other graphs of this kind have been plotted for other values of parameter $\beta$. In all of these graphs, we have observed that $\partial \text{Ei}_{\alpha,\beta}(x)/\partial\beta$ tend to zero for large values of parameter $\alpha$.

For the rest of results given in Tables 1, 2, and 3, we have derived the following reduction formulas.

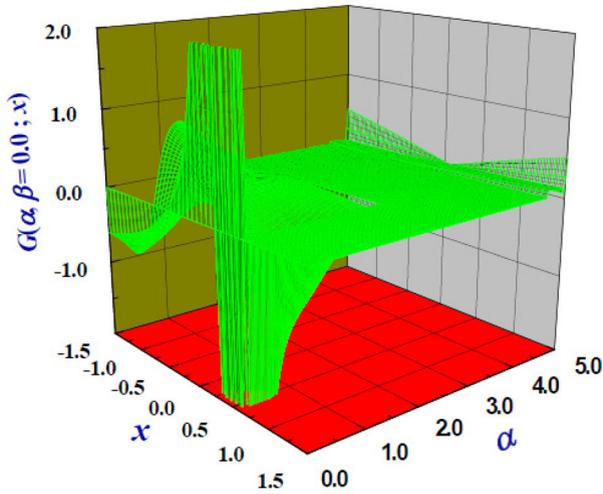

**Fig. 1** Differentiation of the integral Mittag-Leffler function with respect to parameter $\alpha$ as a function of variable $x$ and parameter $\alpha$ at $\beta = 0$

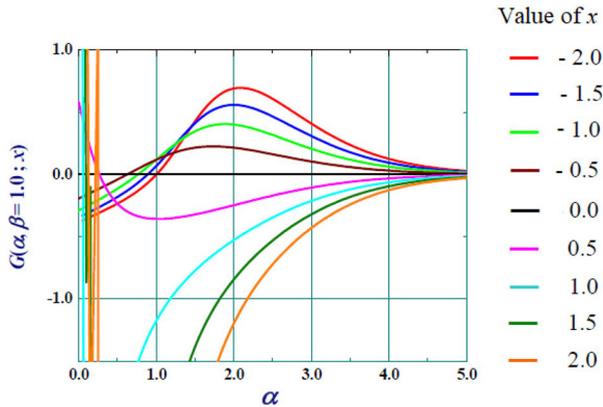

**Fig. 2** Differentiation of the integral Mittag-Leffler function with respect to parameter $\alpha$ as a function of parameter $\alpha$ at $\beta = 1$

**Definition 2.3** According to [4, Eqn. 6.2.1(63)], define the function:

$$Q(a, x) := \sum_{k=1}^{\infty} \frac{x^k}{(a)_k} \psi(k + a)$$
$$= e^x \left[ x^{1-a} \psi(a) \gamma(a, x) + \frac{x}{a^2} {}_2F_2\left( \begin{matrix} a, a \\ a+1, a+1 \end{matrix} \bigg| -x \right) \right]. \tag{2.36}$$

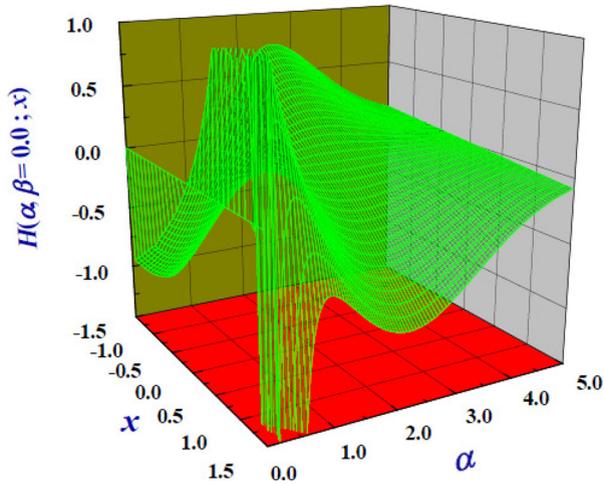

**Fig. 3** Differentiation of the integral Mittag-Leffler function with respect to parameter $\beta$ as a function of variable $x$ and parameter $\alpha$ at $\beta = 0$

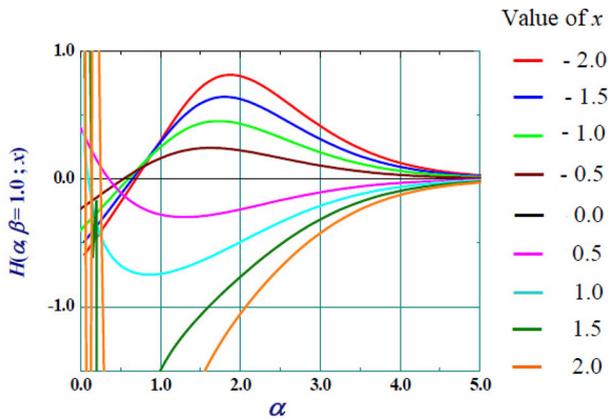

**Fig. 4** Differentiation of the integral Mittag-Leffler function with respect to parameter $\beta$ as a function of parameter $\alpha$ at $\beta = 1$

**Theorem 2.4** *For $q \in \mathbb{N}$ and arbitrary $\beta$, the following reduction formula holds true:*

$$\left.\frac{\partial \mathrm{Ei}_{\alpha,\beta}(x)}{\partial \alpha}\right|_{\alpha=1/q} = -\sum_{h=0}^{q-1} \frac{x^{-h}}{\Gamma\left(\beta - \frac{h}{q}\right)} Q\left(\beta - \frac{h}{q}, x^q\right). \tag{2.37}$$

*Proof* Apply the series decomposition

$$\sum_{k=1}^{\infty} a(k) = \sum_{h=0}^{q-1} \sum_{k=1}^{\infty} a(qk - h), \tag{2.38}$$

**Table 1** First derivative of the integral Mittag-Leffler function $\text{Ei}_{\alpha,\beta}(x)$ with respect to parameter $\alpha$ derived for some values of parameters $\alpha$ and $\beta$

| $\alpha$ | $\beta$ | $\dfrac{\partial \text{Ei}_{\alpha,\beta}(x)}{\partial \alpha}\ (x \neq 0)$ |
|---|---|---|
| 0 | $\beta$ | $\dfrac{x\,\psi(\beta)}{(x-1)\,\Gamma(\beta)}$ |
| $\tfrac{1}{3}$ | 1 | $e^{x^3}\left[-\dfrac{9x\,{}_2F_2\left(\tfrac{1}{3},\tfrac{1}{3};\tfrac{4}{3},\tfrac{4}{3};-x^3\right)+\psi\left(\tfrac{1}{3}\right)\gamma\left(\tfrac{1}{3},x^3\right)}{\Gamma\left(\tfrac{1}{3}\right)} - \dfrac{9x^2\,{}_2F_2\left(\tfrac{2}{3},\tfrac{2}{3};\tfrac{5}{3},\tfrac{5}{3};-x^3\right)+4\psi\left(\tfrac{2}{3}\right)\gamma\left(\tfrac{2}{3},x^3\right)}{4\Gamma\left(\tfrac{2}{3}\right)} + \text{Ei}\left(-x^3\right) - 3\ln x\right] - \gamma$ |
| $\tfrac{1}{3}$ | 2 | $e^{x^3}\left[-\dfrac{9x\,{}_2F_2\left(\tfrac{4}{3},\tfrac{4}{3};\tfrac{7}{3},\tfrac{7}{3};-x^3\right)}{16\,\Gamma\left(\tfrac{4}{3}\right)} - \dfrac{9x^2\,{}_2F_2\left(\tfrac{5}{3},\tfrac{5}{3};\tfrac{8}{3},\tfrac{8}{3};-x^3\right)}{25\,\Gamma\left(\tfrac{5}{3}\right)}\right] - \dfrac{1}{x^3}\left[3\ln x + \gamma - 1 - \text{Ei}(-x^3) - (\gamma-1)\gamma(2,x^3) + \dfrac{\psi\left(\tfrac{5}{3}\right)\gamma\left(\tfrac{5}{3},x^3\right)}{\Gamma\left(\tfrac{5}{3}\right)} + \dfrac{\psi\left(\tfrac{4}{3}\right)\gamma\left(\tfrac{4}{3},x^3\right)}{\Gamma\left(\tfrac{4}{3}\right)}\right] - \dfrac{1}{x^3}$ |
| $\tfrac{1}{3}$ | 3 | $-\dfrac{e^{x^3}}{2x^6}\left[\dfrac{9x^8\,{}_2F_2\left(\tfrac{8}{3},\tfrac{8}{3};\tfrac{11}{3},\tfrac{11}{3};-x^3\right)+64\psi\left(\tfrac{8}{3}\right)\gamma\left(\tfrac{8}{3},x^3\right)}{32\,\Gamma\left(\tfrac{8}{3}\right)} - \dfrac{18x^7\,{}_2F_2\left(\tfrac{7}{3},\tfrac{7}{3};\tfrac{10}{3},\tfrac{10}{3};-x^3\right)+98\psi\left(\tfrac{7}{3}\right)\gamma\left(\tfrac{7}{3},x^3\right)}{49\,\Gamma\left(\tfrac{7}{3}\right)}\right. $ $\left. + e^{-x^3}\left(x^3+3\right) - 2\text{Ei}\left(-x^3\right) + \left(\tfrac{3}{2}-\gamma\right)\gamma\left(3,x^3\right) + 6\ln x + 2\gamma - 3\right]$ |
| $\tfrac{1}{3}$ | $\beta$ | $e^{-x^3}\left[x\,\Gamma\!\left(\beta-\tfrac{2}{3}\right)\,{}_2\mathbf{F}_2\!\left(\beta-\tfrac{2}{3},\beta-\tfrac{2}{3};\beta+\tfrac{1}{3},\beta+\tfrac{1}{3};-x^3\right) + x^2\,\Gamma\!\left(\beta-\tfrac{1}{3}\right)\,{}_2\mathbf{F}_2\!\left(\beta-\tfrac{1}{3},\beta-\tfrac{1}{3};\beta+\tfrac{2}{3},\beta+\tfrac{2}{3};-x^3\right)\right.$ $\left. + \dfrac{\psi\!\left(\beta-\tfrac{2}{3}\right)\gamma\!\left(\beta-\tfrac{2}{3},x^3\right)}{\Gamma\!\left(\beta-\tfrac{2}{3}\right)} + \dfrac{\psi\!\left(\beta-\tfrac{1}{3}\right)\gamma\!\left(\beta-\tfrac{1}{3},x^3\right)}{\Gamma\!\left(\beta-\tfrac{1}{3}\right)} + \dfrac{\psi(\beta)\gamma(\beta,x^3)}{\Gamma(\beta)}\right]$ $+ x^3\,\Gamma(\beta)\,{}_2\mathbf{F}_2\!\left(\beta,\beta;\beta+1,\beta+1;-x^3\right) + x^{3(1-\beta)}$ |
| $\tfrac{1}{2}$ | 1 | $e^{x^2}\left[-\dfrac{4x}{\sqrt{\pi}}\,{}_2F_2\left(\tfrac{1}{2},\tfrac{1}{2};\tfrac{3}{2},\tfrac{3}{2};-x^2\right) + (\gamma+\log 4)\,\text{erf}\,x + \text{Ei}\left(-x^2\right) - 2\ln x\right] - \gamma$ |
| $\tfrac{1}{2}$ | 2 | $-e^{x^2}\left(\dfrac{8x\,{}_2F_2\left(\tfrac{3}{2},\tfrac{3}{2};\tfrac{5}{2},\tfrac{5}{2};-x^2\right)}{9\sqrt{\pi}}\right) + x^{-2}\left[e^{-x^2} - \text{Ei}(-x^2) - (\gamma-1)\gamma\left(2,x^2\right) + \dfrac{2\psi\left(\tfrac{3}{2}\right)\gamma\left(\tfrac{3}{2},x^2\right)}{\sqrt{\pi}} + 2\ln x + \gamma - 1\right]$ |
| $\tfrac{1}{2}$ | 3 | $-\dfrac{e^{x^2}\left\{64x^6\,{}_2F_2\left(\tfrac{5}{2},\tfrac{5}{2};\tfrac{7}{2},\tfrac{7}{2};-x^2\right)+25x\left[\sqrt{\pi}\left[3(2\gamma-3)\left(2-\gamma\left(3,x^2\right)\right)+24\ln x - 12\text{Ei}(-x^2)\right]+16\psi\left(\tfrac{5}{2}\right)\gamma\left(\tfrac{5}{2},x^2\right)\right]+150\sqrt{\pi}x\left(x^2+3\right)\right\}}{300\sqrt{\pi}x^5}$ |
| $\tfrac{1}{2}$ | $\beta$ | $e^{x^2}\left[-\dfrac{x^2}{\Gamma(\beta)}\left(\dfrac{2F_2\left(\beta,\beta;\beta+1,\beta+1;-x^2\right)}{\beta^2} + \psi(\beta)\,x^{-2\beta}\,\gamma\left(\beta,x^2\right)\right)\right.$ $\left. - \dfrac{x^2\,{}_2F_2\left(\beta-\tfrac{1}{2},\beta-\tfrac{1}{2};\beta+\tfrac{1}{2},\beta+\tfrac{1}{2};-x^2\right)}{\left(\beta-\tfrac{1}{2}\right)^2} + \psi\!\left(\beta-\tfrac{1}{2}\right)x^{3-2\beta}\,\gamma\!\left(\beta-\tfrac{1}{2},x^2\right)\right]$ $- \dfrac{1}{x\,\Gamma\!\left(\beta-\tfrac{1}{2}\right)}$ |

**Table 2** First derivative of the integral Mittag-Leffler function $\mathrm{Ei}_{\alpha,\beta}(x)$ with respect to parameter $\alpha$ derived for some values of parameters $\alpha$ and $\beta$

| $\alpha$ | $\beta$ | $\dfrac{\partial \mathrm{Ei}_{\alpha,\beta}(x)}{\partial \alpha}$ $\quad (x \neq 0)$ |
|---|---|---|
| 1 | 0 | $-xe^x \left[\ln x - \mathrm{Ei}(-x)\right]$ |
| 1 | 1 | $-e^x \left[\ln x - \mathrm{Ei}(-x)\right] - \gamma$ |
| 1 | 2 | $1 - \dfrac{\gamma(x+1) + e^x [\ln x - \mathrm{Ei}(-x)]}{x}$ |
| 1 | 3 | $\dfrac{(3x+4)x - 2\gamma[(x+2)x+2] - 4e^x[\ln x - \mathrm{Ei}(-x)]}{4x^2}$ |
| 1 | $\beta$ | $-\dfrac{e^x}{\Gamma(\beta)} \left( \dfrac{{}_2F_2(\beta,\beta;\beta+1,\beta+1;-x)}{\beta^2} + \psi(\beta) \right) x^{-\beta} \gamma(\beta,x)$ |
| 2 | $-1$ | $\dfrac{x}{4} e^{-\sqrt{x}} \left\{ e^{-2\sqrt{x}} \left[2\mathrm{E}_1\left(\sqrt{x}\right) + \ln x\right] + 2\mathrm{Ei}\left(\sqrt{x}\right) - \ln x \right\}$ |
| 2 | 0 | $\sqrt{x} \left\{ \sinh\sqrt{x} \left[\mathrm{Chi}\sqrt{x} - \ln\sqrt{x}\right] - \mathrm{Shi}\sqrt{x}\cosh\sqrt{x} \right\}$ |
| 2 | 1 | $\cosh\sqrt{x} \left[\mathrm{Chi}\sqrt{x} - \ln\sqrt{x}\right] - \mathrm{Shi}\sqrt{x}\sinh\sqrt{x} - \gamma$ |
| 2 | 2 | $\dfrac{1}{\sqrt{x}} \left\{ \sinh\sqrt{x} \left[\mathrm{Chi}\sqrt{x} - \ln\sqrt{x}\right] - \mathrm{Shi}\sqrt{x}\cosh\sqrt{x} \right\} + 1 - \gamma$ |
| 3 | $-2$ | $-\dfrac{x}{9}\left\{ e^{\sqrt[3]{x}} \left[3\mathrm{E}_1\left(\sqrt[3]{x}\right) + \ln x\right] + e^{-\sqrt[3]{-1}\sqrt[3]{x}} \left[3\mathrm{E}_1\left(-\sqrt[3]{-1}\sqrt[3]{x}\right) + \ln x - 2\pi i\right] \right.$ $\left. + e^{(-1)^{2/3}\sqrt[3]{x}} \left[3\mathrm{E}_1\left((-1)^{2/3}\sqrt[3]{x}\right) + \ln x + 2\pi i\right] \right\}$ |
| 3 | $-1$ | $-\dfrac{x^{2/3}}{3} \left\{ e^{\sqrt[3]{x}} \left[\mathrm{E}_1\left(\sqrt[3]{x}\right) + \dfrac{\ln x}{3}\right] - \dfrac{1}{3}\sqrt[3]{-1}e^{(-1)^{2/3}\sqrt[3]{x}} \left[3\mathrm{E}_1\left((-1)^{2/3}\sqrt[3]{x}\right) + \ln x + 2\pi i\right] \right.$ $\left. + \dfrac{1}{3}(-1)^{2/3} e^{-\sqrt[3]{-1}\sqrt[3]{x}} \left[3\Gamma\left(0,-\sqrt[3]{-1}\sqrt[3]{x}\right) + \ln x - 2\pi i\right] \right\}$ |
| 3 | 0 | $-\dfrac{\sqrt[3]{x}}{9} \left\{ (-1)^{5/6} e^{-\sqrt[3]{-1}\sqrt[3]{x}} \left[3i\mathrm{E}_1\left(-\sqrt[3]{-1}\sqrt[3]{x}\right) + i\ln x + 2\pi\right] + e^{\sqrt[3]{x}} \left[3\mathrm{E}_1\left(\sqrt[3]{x}\right) + \ln x\right] \right.$ $\left. + (-1)^{2/3} e^{(-1)^{2/3}\sqrt[3]{x}} \left[3\mathrm{E}_1\left((-1)^{2/3}\sqrt[3]{x}\right) + \ln x + 2\pi i\right] \right\}$ |
| 4 | $-2$ | $\dfrac{x^{3/4}}{2} \left\{ \mathrm{Chi}\sqrt[4]{x}\sinh\sqrt[4]{x} + \mathrm{Ci}\sqrt[4]{x}\sin\sqrt[4]{x} - \mathrm{Shi}\sqrt[4]{x}\cosh\sqrt[4]{x} - \mathrm{Si}\sqrt[4]{x}\cos\sqrt[4]{x} - \ln\sqrt[4]{x}\left[\sin\sqrt[4]{x} + \sinh\sqrt[4]{x}\right] \right\}$ |
| 4 | 0 | $\dfrac{1}{2}\sqrt[4]{x} \left\{ \mathrm{Chi}\sqrt[4]{x}\sinh\sqrt[4]{x} + \sin\sqrt[4]{x}\left[\ln\sqrt[4]{x} - \mathrm{Ci}\sqrt[4]{x}\right] - \mathrm{Shi}\sqrt[4]{x}\cosh\sqrt[4]{x} + \mathrm{Si}\sqrt[4]{x}\cos\sqrt[4]{x} - \ln\sqrt[4]{x}\sinh\sqrt[4]{x} \right\}$ |

**Table 3** First derivative of the integral Mittag-Leffler function $\mathrm{Ei}_{\alpha,\beta}(x)$ with respect to parameter $\beta$ derived for some values of parameters $\alpha$ and $\beta$

| $\alpha$ | $\beta$ | $\dfrac{\partial \mathrm{Ei}_{\alpha,\beta}(x)}{\partial \beta}$  $(x \neq 0)$ |
|---|---|---|
| $\frac{1}{3}$ | 0 | $\dfrac{1}{3}\left[-\dfrac{e^{x^3}\left[9x\,_2F_2\left(\frac{1}{3},\frac{1}{3};\frac{4}{3},\frac{4}{3};-x^3\right)+\psi\left(\frac{1}{3}\right)\gamma\left(\frac{1}{3},x^3\right)\right]-9x\,_2F_2\left(\frac{1}{3},1;\frac{4}{3},\frac{4}{3};x^3\right)}{\Gamma\left(\frac{1}{3}\right)} - \dfrac{e^{x^3}\left[9x^2\,_2F_2\left(\frac{2}{3},\frac{2}{3};\frac{5}{3},\frac{5}{3};-x^3\right)+4\psi\left(\frac{2}{3}\right)\gamma\left(\frac{2}{3},x^3\right)\right]-9x^2\,_2F_2\left(\frac{2}{3},1;\frac{5}{3},\frac{5}{3};x^3\right)}{4\Gamma\left(\frac{2}{3}\right)}\right.$ $\left.+\mathrm{Ei}(x^3)+e^{x^3}\left[\mathrm{Ei}(-x^3)-3\ln x\right]-3\ln x-2\gamma\right]$ |
| $\frac{1}{2}$ | 0 | $\dfrac{2x}{\sqrt{\pi}}\left\{2F_2\left(\frac{1}{2},1;\frac{3}{2},\frac{3}{2};x^2\right)-e^{x^2}\,_2F_2\left(\frac{1}{2},\frac{1}{2};\frac{3}{2},\frac{3}{2};-x^2\right)\right\}$ $+\dfrac{1}{2}\left\{e^{x^2}\left[\mathrm{Chi}\left(x^2\right)+(\gamma+\ln 4)\,\mathrm{erf}\,x-\mathrm{Shi}\left(x^2\right)-2\ln x\right]+\mathrm{Chi}\left(x^2\right)+\mathrm{Shi}\left(x^2\right)-2(\gamma+\ln x)\right\}$ |
| 1 | 0 | $e^x\,\mathrm{Ei}(-x)+\mathrm{Ei}(x)-(e^x+1)\ln x-2\gamma$ |
| 2 | 0 | $2\,\mathrm{Chi}\sqrt{x}\,[\cosh\sqrt{x}+1]-2\,\mathrm{Shi}\sqrt{x}\,\sinh\sqrt{x}-\ln x\,[\cosh\sqrt{x}+1]-4\gamma$ |
| 2 | 1 | $\dfrac{x}{4}\,_3F_4\left(1,1,1;2,2,2,\frac{3}{2};\frac{x}{4}\right)+\left[\ln\sqrt{x}-\mathrm{Chi}\sqrt{x}\right]^2-\mathrm{Shi}^2\sqrt{x}-\gamma^2$ |
| 4 | 0 | $2\,\mathrm{Chi}\sqrt[4]{x}\,[\cos\sqrt[4]{x}+1]+2\,\mathrm{Ci}\sqrt[4]{x}\,[\cos\sqrt[4]{x}+1]-2\,\mathrm{Shi}\sqrt[4]{x}\,\sinh\sqrt[4]{x}+2\,\mathrm{Si}\sqrt[4]{x}\,\sin\sqrt[4]{x}-\dfrac{1}{2}\ln x\,[\cos\sqrt[4]{x}+\cosh\sqrt[4]{x}+2]-8\gamma$ |

to obtain

$$\frac{\partial \mathrm{Ei}_{\alpha,\beta}(x)}{\partial \alpha} = -\sum_{h=0}^{q-1}\sum_{k=1}^{\infty} \frac{x^{qk-h}\psi(\alpha(qk-h)+\beta)}{\Gamma(\alpha(qk-h)+\beta)}. \quad (2.39)$$

Therefore, taking into account (2.36), we arrive at

$$\left.\frac{\partial \mathrm{Ei}_{\alpha,\beta}(x)}{\partial \alpha}\right|_{\alpha=1/q} = -\sum_{h=0}^{q-1} \frac{x^{-h}}{\Gamma\left(\beta-\frac{h}{q}\right)} \sum_{k=1}^{\infty} \frac{(x^q)^k \psi\left(k+\beta-\frac{h}{q}\right)}{\left(\beta-\frac{h}{q}\right)_k}, \quad (2.40)$$

which is equivalent to (2.37), as we wanted to prove. □

**Definition 2.5** According to the definition of the generalized hypergeometric function, and the property

$$(k+a)(a)_{k+1} = \frac{a^2}{(a)_k}\left[(a+1)_k\right]^2, \quad (2.41)$$

define the function

$$S(a,x) := \sum_{k=0}^{\infty} \frac{x^{k+1}}{(k+a)(a)_{k+1}}$$
$$= \frac{x}{a^2} \sum_{k=0}^{\infty} \frac{x^k (a)_k (1)_k}{k!\left[(a+1)_k\right]^2}$$
$$= \frac{x}{a^2}\, {}_2F_2\left(\begin{array}{c}1, a\\1+a, 1+a\end{array}\middle| x\right). \quad (2.42)$$

**Theorem 2.6** *For $q \in \mathbb{N}$, the following reduction formula holds true:*

$$\left.\frac{\partial \mathrm{Ei}_{\alpha,\beta}(x)}{\partial \beta}\right|_{\alpha=1/q,\,\beta=0} = \frac{1}{q}\sum_{h=0}^{q-1} \frac{x^{-h}}{\Gamma\left(1-\frac{h}{q}\right)}\left[S\left(1-\frac{h}{q},x^q\right) - Q\left(1-\frac{h}{q},x^q\right)\right]. \quad (2.43)$$

*Proof* From (2.5) and taking into account (2.3), we have

$$\left.\frac{\partial \mathrm{Ei}_{\alpha,\beta}(x)}{\partial \beta}\right|_{\beta=0} = -\alpha \sum_{k=1}^{\infty} \frac{x^k \psi(\alpha k)}{\alpha k\, \Gamma(\alpha k)}$$
$$= -\alpha \sum_{k=1}^{\infty} \frac{x^k \psi(\alpha k+1)}{\Gamma(\alpha k+1)} + \sum_{k=1}^{\infty} \frac{x^k}{k\,\Gamma(\alpha k+1)}. \quad (2.44)$$

Apply now (2.38),

$$\left.\frac{\partial \mathrm{Ei}_{\alpha,\beta}(x)}{\partial \beta}\right|_{\beta=0} = -\alpha \sum_{h=0}^{q-1} \sum_{k=1}^{\infty} \frac{x^{qk-h} \psi(\alpha(qk-h)+1)}{\Gamma(\alpha(qk-h)+1)}$$
$$+ \sum_{h=0}^{q-1} \sum_{k=1}^{\infty} \frac{x^{qk-h}}{(qk-h) \Gamma(\alpha(qk-h)+1)}, \quad (2.45)$$

thereby, taking $\alpha = 1/q$ with $q \in \mathbb{N}$ and recalling (2.36), we arrive at

$$\left.\frac{\partial \mathrm{Ei}_{\alpha,\beta}(x)}{\partial \beta}\right|_{\alpha=1/q,\,\beta=0} = -\frac{1}{q} \sum_{h=0}^{q-1} x^{-h} \sum_{k=1}^{\infty} \frac{(x^q)^k \psi\left(k - \frac{h}{q} + 1\right)}{\Gamma\left(k - \frac{h}{q} + 1\right)}$$
$$+ \frac{1}{q} \sum_{h=0}^{q-1} x^{-h} \sum_{k=1}^{\infty} \frac{(x^q)^k}{\left(k - \frac{h}{q}\right) \Gamma\left(k - \frac{h}{q} + 1\right)}$$
$$= \frac{1}{q} \sum_{h=0}^{q-1} \frac{x^{-h}}{\Gamma\left(1 - \frac{h}{q}\right)} \left[ -Q\left(1 - \frac{h}{q}, x^q\right) + \sum_{k=0}^{\infty} \frac{(x^q)^{k+1}}{\left(k + 1 - \frac{h}{q}\right) \left(1 - \frac{h}{q}\right)_{k+1}} \right], \quad (2.46)$$

which, according to (2.42), is equivalent to (2.6). □

In Table 1, note that $_p\mathbf{F}_q\left(a_1, \ldots, a_p; b_1 \ldots, b_q; z\right)$ denotes the function [16, Eqn. 16.2.5]:

$$_p\mathbf{F}_q\left(\begin{array}{c}a_1,\ldots,a_p\\b_1,\ldots,b_q\end{array}\bigg| z\right) = \sum_{k=0}^{\infty} \frac{(a_1)\cdots(a_p)}{\Gamma(b_1+k)\cdots\Gamma(b_q+k)} \frac{z^k}{k!}. \quad (2.47)$$

## 3 Differentiation of the integral Wright function with respect to parameters

Considered as a kind of generalization of the Bessel functions at the beginning, today the Wright function (called also Bessel-Maitland function), which was introduced in 1933 and 1940 [20, 21],

$$\mathrm{W}_{\alpha,\beta}(x) = \sum_{k=0}^{\infty} \frac{x^k}{k!\,\Gamma(\alpha k + \beta)}, \quad (3.1)$$

is an important special function involved in solutions of various physical problems such as space diffusion models, stochastic processes, probability distributions and other natural phenomena [10, 13]. The infinite series in (3.1) differs from the Mittag-Leffler function by appearance of factorial only. Therefore, similarly to (2.1) and (2.5)

but using (1.4), we have

$$I(\alpha, \beta; x) = \frac{\partial \mathrm{Wi}_{\alpha,\beta}(x)}{\partial \alpha} = -\sum_{k=1}^{\infty} \frac{\psi(\alpha k + \beta)}{k!\, \Gamma(\alpha k + \beta)} x^k, \qquad (3.2)$$

$$J(\alpha, \beta; x) = \frac{\partial \mathrm{Wi}_{\alpha,\beta}(x)}{\partial \beta} = -\sum_{k=1}^{\infty} \frac{\psi(\alpha k + \beta)}{k\, k!\, \Gamma(\alpha k + \beta)} x^k. \qquad (3.3)$$

From (3.1), it follows that

$$\frac{\partial W_{\alpha,\beta}(x)}{\partial \beta} = -\sum_{k=0}^{\infty} \frac{\psi(\alpha k + \beta)}{k!\, \Gamma(\alpha k + \beta)} x^k, \qquad (3.4)$$

and comparing (3.3) to (3.4), first derivatives with respect to parameter $\beta$ of the Wright function and the integral Wright functions are interrelated by

$$\frac{\partial \mathrm{Wi}_{\alpha,\beta}(x)}{\partial \alpha} = \begin{cases} \frac{\psi(\beta)}{\Gamma(\beta)} + \frac{\partial W_{\alpha,\beta}(x)}{\partial \beta}, & \beta \neq 0, \\ \frac{\partial W_{\alpha,\beta}(x)}{\partial \beta}, & \beta = 0. \end{cases} \qquad (3.5)$$

Also, comparing (3.2) to (3.3), we obtain

$$\frac{\partial \mathrm{Wi}_{\alpha,\beta}(x)}{\partial \alpha} = x \frac{\partial}{\partial x} \left( \frac{\partial \mathrm{Wi}_{\alpha,\beta}(x)}{\partial \beta} \right), \qquad (3.6)$$

or, in equivalent form,

$$\frac{\partial \mathrm{Wi}_{\alpha,\beta}(x)}{\partial \beta} = \int_0^x \left( \frac{\partial \mathrm{Wi}_{\alpha,\beta}(t)}{\partial \alpha} \right) \frac{dt}{t}. \qquad (3.7)$$

Note that $I(\alpha, \beta; 0) = J(\alpha, \beta; 0) = 0$, so that we will consider $x \neq 0$ in this section in order to avoid indeterminate expressions.

**Theorem 3.1** *For $\alpha = 1$ and $\beta \notin \mathbb{Z}$, the following reduction formula holds true:*

$$\begin{aligned}
\left. \frac{\partial \mathrm{Wi}_{\alpha,\beta}(x)}{\partial \alpha} \right|_{\alpha=1} &= \frac{\psi(\beta)}{\Gamma(\beta)} - x^{(1-\beta)/2} I_{\beta-1}(2\sqrt{x}) \left[ \psi(\beta) + \frac{x}{\beta(\beta-1)}\, {}_3F_4\left( \begin{array}{c} 1, 1, 3/2 \\ 2, 2, 1+\beta, 1+\beta \end{array} \bigg| 4x \right) \right] \\
&\quad - \frac{\Gamma(1-\beta)}{\beta\, \Gamma(1+\beta)} x^{(1+\beta)/2} I_{1-\beta}(2\sqrt{x})\, {}_2F_3\left( \begin{array}{c} \beta, \frac{1}{2}+\beta \\ 2\beta, 1+\beta, 1+\beta \end{array} \bigg| 4x \right).
\end{aligned} \qquad (3.8)$$

**Proof** From the following result for $\beta \notin \mathbb{Z}$ [14]:

$$\sum_{k=1}^{\infty} [\psi(k+\beta) - \psi(\beta)] \frac{x^k}{(\beta)_k \, k!}$$
$$= \frac{x}{\beta(\beta-1)} \left[ {}_0F_1\left(\left.\begin{array}{c}-\\ \beta\end{array}\right| x\right) {}_3F_4\left(\left.\begin{array}{c}1, 1, 3/2\\ 2, 2, 1+\beta, 1+\beta\end{array}\right| 4x\right) \right.$$
$$\left. - \frac{1}{\beta} {}_0F_1\left(\left.\begin{array}{c}-\\ 2-\beta\end{array}\right| x\right) {}_2F_3\left(\left.\begin{array}{c}\beta, \frac{1}{2}+\beta\\ 2\beta, 1+\beta, 1+\beta\end{array}\right| 4x\right) \right],$$

and taking into account [12, Sect. 9.14]

$${}_0F_1\left(\left.\begin{array}{c}-\\ \beta\end{array}\right| x\right) = \Gamma(\beta) \, x^{(1-\beta)/2} I_{\beta-1}\left(2\sqrt{x}\right),$$

after simplification, we arrive at (3.8), as we wanted to prove. □

Next, we present some reduction formulas of the derivative with respect to $\alpha$ of the integral Wright function for $\alpha = 1$ and $\beta = 0, 1, 2$.

**Theorem 3.2** *For $\alpha = 1$ and $\beta = 0$, the following reduction formula holds true:*

$$\left.\frac{\partial \mathrm{Wi}_{\alpha,\beta}(x)}{\partial \alpha}\right|_{\alpha=1,\beta=0}$$
$$= \frac{1}{2}\left\{I_0\left(2\sqrt{x}\right) - \sqrt{x}\left[\ln x \, I_1\left(2\sqrt{x}\right) - 2K_1\left(2\sqrt{x}\right)\right]\right\} - 1. \quad (3.9)$$

**Proof** Apply (3.5), taking into account (6.6) (see Appendix 6). □

**Theorem 3.3** *For $\alpha = 1$ and $\beta = 1$, the following reduction formula holds true:*

$$\left.\frac{\partial \mathrm{Wi}_{\alpha,\beta}(x)}{\partial \alpha}\right|_{\alpha=\beta=1} = -\gamma - \frac{1}{2}\ln x \, I_0\left(2\sqrt{x}\right) - K_0\left(2\sqrt{x}\right). \quad (3.10)$$

**Proof** Apply (3.5), taking into account (6.8) and knowing that $\psi(1) = -\gamma$. □

**Theorem 3.4** *For $\alpha = 1$ and $\beta = 2$, the following reduction formula holds true:*

$$\left.\frac{\partial \mathrm{Wi}_{\alpha,\beta}(x)}{\partial \alpha}\right|_{\alpha=1, \beta=2} = 1 - \gamma - \frac{I_0\left(2\sqrt{x}\right)}{2x} + \frac{2K_1\left(2\sqrt{x}\right) - \ln x \, I_1\left(2\sqrt{x}\right)}{2\sqrt{x}}. \quad (3.11)$$

**Proof** Apply (3.5), taking into account (6.10) and knowing that $\psi(2) = 1 - \gamma$. □

Finally, we present some reduction formulas of the derivative with respect to $\beta$ of the integral Wright function for $\alpha = 1$ and $\beta = 0, 1, 2$.

**Theorem 3.5** *For $\alpha = 1$ and $\beta = 0$, the following reduction formula holds true:*

$$\left.\frac{\partial \mathrm{Wi}_{\alpha,\beta}(x)}{\partial \beta}\right|_{\alpha=1,\,\beta=0} = -\gamma - \frac{1}{2}\ln x\, I_0\left(2\sqrt{x}\right) - K_0\left(2\sqrt{x}\right) + x\,{}_2F_3\left(\begin{array}{c}1,1\\2,2,2\end{array}\bigg|\,x\right).$$
(3.12)

*Proof* From (3.3) and (2.3) we have

$$\left.\frac{\partial \mathrm{Wi}_{\alpha,\beta}(x)}{\partial \beta}\right|_{\alpha=1,\,\beta=0} = -\sum_{k=1}^{\infty}\frac{x^k \psi(k)}{k\,k!\,\Gamma(k)} = -\sum_{k=1}^{\infty}\frac{x^k}{(k!)^2}\left[\psi(k+1) - \frac{1}{k}\right]$$

$$= -\gamma - \sum_{k=0}^{\infty}\frac{x^k}{(k!)^2}\psi(k+1) + x\sum_{k=0}^{\infty}\frac{x^k}{(k+1)\left[(k+1)!\right]^2}.$$
(3.13)

However,

$$\sum_{k=0}^{\infty}\frac{x^k}{(k+1)\left[(k+1)!\right]^2} = \sum_{k=0}^{\infty}\frac{x^k\left[(1)_k\right]^2}{k!\left[(2)_k\right]^3} = {}_2F_3\left(\begin{array}{c}1,1\\2,2,2\end{array}\bigg|\,x\right),$$
(3.14)

so that, from (6.4) and (3.14), we finally obtain (3.12), as we wanted to prove. □

**Theorem 3.6** *For $\alpha = 1$ and $\beta = 1$, the following reduction formula holds true:*

$$\left.\frac{\partial \mathrm{Wi}_{\alpha,\beta}(x)}{\partial \beta}\right|_{\alpha=\beta=1}$$
$$= \frac{x}{2}\left[{}_3F_4\left(\begin{array}{c}1,1,1\\2,2,2,2\end{array}\bigg|\,x\right) - \ln x\,{}_2F_3\left(\begin{array}{c}1,1\\2,2,2\end{array}\bigg|\,x\right)\right] + \frac{1}{2}G_{1,3}^{3,0}\left(x,\bigg|\begin{array}{c}1\\0,0,0\end{array}\right)$$
$$- \frac{1}{4}\ln x\,(\ln x + 4\gamma) - \gamma^2 - \frac{\pi^2}{12},$$
(3.15)

*where $G_{p,q}^{m,n}\left(x\,\bigg|\begin{array}{c}a_1,\ldots,a_n,a_{n+1},\ldots,a_p\\b_1,\ldots,b_m,b_{m+1},\ldots,b_q\end{array}\right)$ denotes the Meijer-G function [16, Sect. 16.17].*

*Proof* Apply (3.7), taking into account (3.10), to obtain

$$\left.\frac{\partial \mathrm{Wi}_{\alpha,\beta}(x)}{\partial \beta}\right|_{\alpha=\beta=1} = \int_0^x \left.\frac{\partial \mathrm{Wi}_{\alpha,\beta}(t)}{\partial \alpha}\right|_{\alpha=\beta=1}\frac{dt}{t}$$
$$= \int_0^x \left(-\gamma - \frac{1}{2}\ln t\,I_0\left(2\sqrt{t}\right) - K_0\left(2\sqrt{t}\right)\right)\frac{dt}{t}.$$
(3.16)

Calculate the above integral with the aid of MATHEMATICA to arrive at (3.15). □

**Theorem 3.7** *For $\alpha = 1$ and $\beta = 2$, the following reduction formula holds true:*

$$\left.\frac{\partial \text{Wi}_{\alpha,\beta}(x)}{\partial \beta}\right|_{\alpha=1,\beta=2}$$
$$= \frac{x}{4}\left[{}_3F_4\left(\begin{array}{c}1,1,1\\2,2,2,2\end{array}\bigg|x\right) - \ln x \, {}_2F_3\left(\begin{array}{c}1,1\\2,2,2\end{array}\bigg|x\right) - \frac{1}{2}{}_2F_3\left(\begin{array}{c}1,1\\2,3,3\end{array}\bigg|x\right)\right]$$
$$-\frac{1}{2}G_{1,3}^{3,0}\left(x,\bigg|\begin{array}{c}1\\-1,0,0\end{array}\right) + \frac{1}{2x} + \frac{1}{2}\ln x\left(1 - 2\gamma - \frac{1}{2}\ln x\right)$$
$$-\gamma^2 - \frac{\pi^2}{12} - \frac{1}{2}. \qquad (3.17)$$

*Proof* Apply (3.7), taking into account (3.17), to obtain

$$\left.\frac{\partial \text{Wi}_{\alpha,\beta}(x)}{\partial \beta}\right|_{\alpha=1,\beta=2} = \int_0^x \left.\frac{\partial \text{Wi}_{\alpha,\beta}(t)}{\partial \alpha}\right|_{\alpha=1,\beta=2} \frac{dt}{t}$$
$$= \int_0^x \left(1 - \gamma - \frac{I_0(2\sqrt{t})}{2t} + \frac{2K_1(2\sqrt{t}) - \ln t \, I_1(2\sqrt{t})}{2\sqrt{t}}\right) \frac{dt}{t}. \qquad (3.18)$$

Calculate the above integral with the aid of MATHEMATICA to arrive at (3.17). □

We summarize the results of this section for the first derivatives of the integral Wright function with respect to $\alpha$ and $\beta$ in Tables 4 and 5.

Figures 5 and 6 present some graphs of the first derivative of the integral Wright function with respect to the first parameter. In Fig. 5, $\partial \text{Wi}_{\alpha,\beta}(x)/\partial\alpha$ is plotted as a function of variable $x$ and parameter $\alpha$ at fixed value of parameter $\beta = 0$. In Fig. 6, $\partial \text{Wi}_{\alpha,\beta}(x)/\partial\alpha$ is plotted as a function of parameter $\alpha$ at fixed values of $x$ and $\beta = 1$. Many other graphs of this kind have been plotted for other values of parameter $\beta$. In all of these graphs, we have observed that $\partial \text{Wi}_{\alpha,\beta}(x)/\partial\alpha$ tend to zero for large values of parameter $\alpha$, as in the case of the integral Mittag-Leffler function.

**Table 4** Derivative of the integral Wright function with respect to the first parameter

| $\alpha$ | $\beta$ | $\frac{\partial \text{Wi}_{\alpha,\beta}(x)}{\partial \alpha}$ ($x \neq 0$) |
|---|---|---|
| 1 | 0 | $-1 + \frac{1}{2}I_0(2\sqrt{x}) + \sqrt{x}\left[K_1(2\sqrt{x}) - \frac{1}{2}\ln x \, I_1(2\sqrt{x})\right]$ |
| 1 | $\frac{1}{2}$ | $\frac{1}{\sqrt{\pi}}\left\{\cosh(2\sqrt{x})\left[\text{Chi}(4\sqrt{x}) - \frac{1}{2}\ln x\right] - \text{Shi}(4\sqrt{x})\sinh(2\sqrt{x}) + \psi\left(\frac{1}{2}\right)\right\}$ |
| 1 | 1 | $-\gamma - \frac{1}{2}\ln x \, I_0(2\sqrt{x}) - K_0(2\sqrt{x})$ |
| 1 | $\frac{3}{2}$ | $\frac{1}{2\sqrt{\pi x}}\left\{\sinh(2\sqrt{x})\left[2\,\text{Chi}(4\sqrt{x}) - \ln x\right] - 2\,\text{Shi}(4\sqrt{x})\cosh(2\sqrt{x}) + 4\sqrt{x}\,\psi\left(\frac{3}{2}\right)\right\}$ |
| 1 | 2 | $1 - \gamma - \frac{I_0(2\sqrt{x})}{2x} + \frac{2K_1(2\sqrt{x}) - \ln x \, I_1(2\sqrt{x})}{2\sqrt{x}}$ |

**Table 5** Derivative of the integral Wright function with respect to the second parameter

| $\alpha$ | $\beta$ | $\dfrac{\partial \mathrm{Wi}_{\alpha,\beta}(x)}{\partial \beta}$ $(x \neq 0)$ |
|---|---|---|
| 1 | 0 | $-\gamma - \tfrac{1}{2}\ln x\, I_0\left(2\sqrt{x}\right) - K_0\left(2\sqrt{x}\right) + x\,{}_2F_3\left(\begin{array}{c}1,1\\2,2,2\end{array}\Big|\,x\right)$ |
| 1 | 1 | $\tfrac{x}{2}\left[{}_3F_4\left(\begin{array}{c}1,1,1\\2,2,2,2\end{array}\Big|\,x\right) - \ln x\,{}_2F_3\left(\begin{array}{c}1,1\\2,2,2\end{array}\Big|\,x\right)\right]$ $+ \tfrac{1}{2}G^{3,0}_{1,3}\left(x,\left|\begin{array}{c}1\\0,0,0\end{array}\right.\right) - \tfrac{1}{4}\ln x\,(\ln x + 4\gamma) - \gamma^2 - \tfrac{\pi^2}{12}$ |
| 1 | 2 | $\tfrac{x}{4}\left[{}_3F_4\left(\begin{array}{c}1,1,1\\2,2,2,2\end{array}\Big|\,x\right) - \ln x\,{}_2F_3\left(\begin{array}{c}1,1\\2,2,2\end{array}\Big|\,x\right) - \tfrac{1}{2}\,{}_2F_3\left(\begin{array}{c}1,1\\2,3,3\end{array}\Big|\,x\right)\right]$ $-\tfrac{1}{2}G^{3,0}_{1,3}\left(x,\left|\begin{array}{c}1\\-1,0,0\end{array}\right.\right) + \tfrac{1}{2x} + \tfrac{1}{2}\ln x\left(1 - 2\gamma - \tfrac{1}{2}\ln x\right) - \gamma^2 - \tfrac{\pi^2}{12} - \tfrac{1}{2}$ |

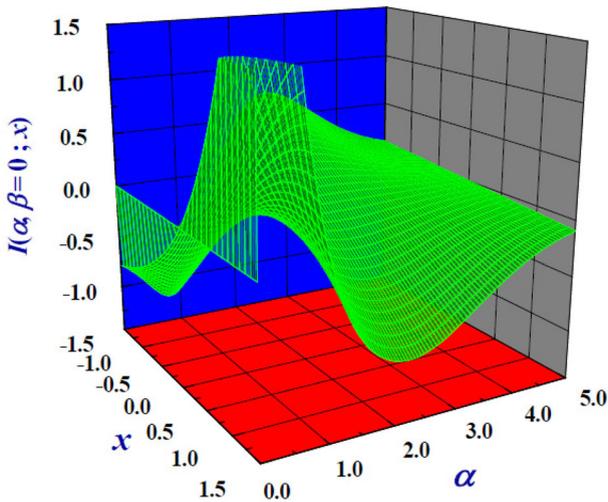

**Fig. 5** Differentiation of the integral Wright function with respect to parameter $\alpha$ as a function of variable $x$ and parameter $\alpha$ at $\beta = 0$

Figures 7 and 8 present some graphs of the first derivative of the integral Wright function with respect to the second parameter. In Fig. 7, $\partial \mathrm{Wi}_{\alpha,\beta}(x)/\partial\beta$ is plotted as a function of variable $x$ and parameter $\alpha$ at fixed value of parameter $\beta = 0$. In Fig. 8, $\partial \mathrm{Wi}_{\alpha,\beta}(x)/\partial\beta$ is plotted as a function of parameter $\alpha$ at fixed values of $x$ and $\beta = 1$. Many other graphs of this kind have been plotted for other values of parameter $\beta$. In all of these graphs, we have observed that $\partial \mathrm{Wi}_{\alpha,\beta}(x)/\partial\beta$ tend to zero for large values of parameter $\alpha$, as in the case of the integral Mittag-Leffler function.

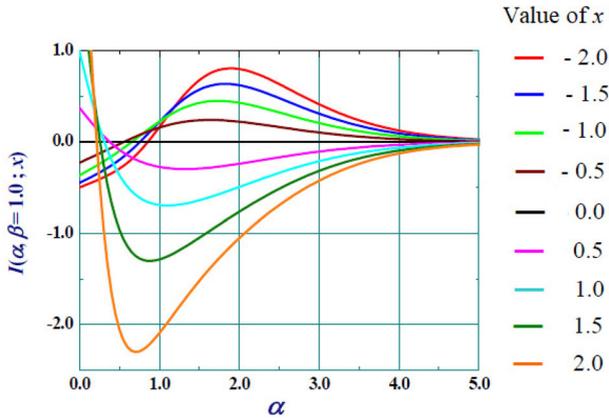

**Fig. 6** Differentiation of the integral Wright function with respect to parameter $\alpha$ as a function of parameter $\alpha$ at $\beta = 1$

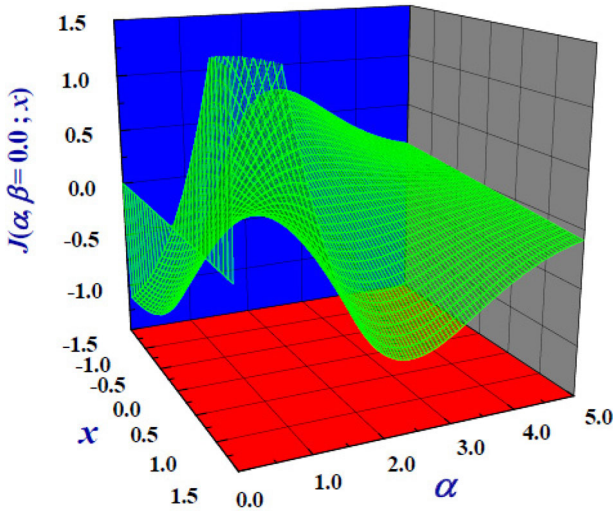

**Fig. 7** Differentiation of the integral Wright function with respect to parameter $\beta$ as a function of variable $x$ and parameter $\alpha$ at $\beta = 0$

## 4 Conclusions

The derivatives with respect to the parameters for the integral Mittag-Leffler function and the integral Wright function have been considered. The differentiation with respect to the parameters of these integral functions has been expressed as infinite sums of quotients of the digamma and gamma functions. On the one hand, these infinite sums have been calculated in closed-form with the help of MATHEMATICA program for some particular cases. On the other hand, some reduction formulas for the derivatives with respect to the parameters have been derived explicitly in order to verify some of the results provided by MATHEMATICA, as well as to provide many other new results. Furthermore, we have represented these infinite sums graphically for particular values of the parameters. Finally, new results are reported for the derivatives with respect to the parameters of the Mittag-Leffler and the Wright functions in the Appendices. In a future work, we will consider the derivatives of the Whittaker and the integral Whittaker functions with respect to the parameters.

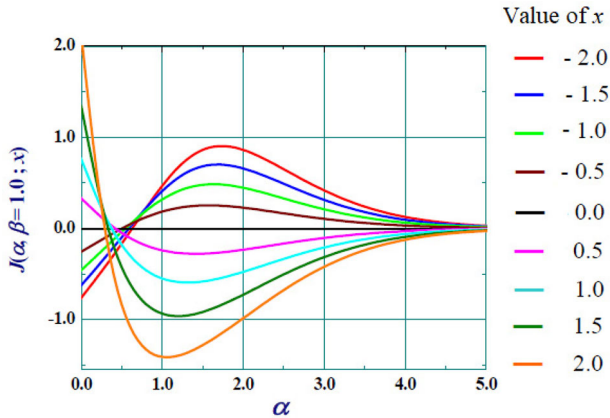

**Fig. 8** Differentiation of the integral Wright function with respect to parameter $\beta$ as a function of parameter $\alpha$ at $\beta = 1$

## 5 Differentiation of two-parameter Mittag-Leffler function with respect to parameters $\alpha$ and $\beta$

As already mentioned in the text (see [2]), derivatives of the Mittag-Leffler function (1.1) with respect to $\alpha$ and $\beta$, i.e.

$$\frac{\partial E_{\alpha,\beta}(x)}{\partial \alpha} = -\sum_{k=0}^{\infty} \frac{\psi(\alpha k + \beta)}{\Gamma(\alpha k + \beta)} k x^k, \tag{5.1}$$

$$\frac{\partial E_{\alpha,\beta}(x)}{\partial \beta} = -\sum_{k=0}^{\infty} \frac{\psi(\alpha k + \beta)}{\Gamma(\alpha k + \beta)} x^k, \tag{5.2}$$

where,

$$\frac{\partial E_{\alpha,\beta}(x)}{\partial \alpha} = x \frac{\partial}{\partial x} \left( \frac{\partial E_{\alpha,\beta}(x)}{\partial \beta} \right), \tag{5.3}$$

can be presented as convolution integrals if the Laplace transform method or the integral representations are used. However, it is possible to present them also in closed-form for integer values of parameters $\alpha$ and $\beta$ [1]. Next, we present some new reduction formulas in closed-form for particular values of $\alpha$ and $\beta$.

**Theorem 5.1** *For $q \in \mathbb{N}$ and arbitrary $\beta$, the following reduction formula holds true:*

$$\left.\frac{\partial E_{\alpha,\beta}(x)}{\partial \beta}\right|_{\alpha=1/q} = -\sum_{h=0}^{q-1} \frac{x^h}{\Gamma\left(\frac{h}{q}+\beta\right)}\left[\psi\left(\frac{h}{q}+\beta\right) + Q\left(\frac{h}{q}+\beta, x^q\right)\right], \quad (5.4)$$

*where recall that the $Q(a,x)$ function is defined in (2.36).*

**Proof** Apply the series decomposition

$$\sum_{k=0}^{\infty} a(k) = \sum_{h=0}^{q-1}\sum_{k=0}^{\infty} a(qk+h), \quad (5.5)$$

to obtain

$$\left.\frac{\partial E_{\alpha,\beta}(x)}{\partial \beta}\right|_{\alpha=1/q} = -\sum_{h=0}^{q-1}\sum_{k=0}^{\infty} \frac{x^{q+hk}\psi(\alpha(qk+h)+\beta)}{\Gamma(\alpha(qk+h)+\beta)}\bigg|_{\alpha=1/q}$$

$$= -\sum_{h=0}^{q-1} \frac{x^h}{\Gamma\left(\frac{h}{q}+\beta\right)} \sum_{k=0}^{\infty} \frac{(x^q)^k \psi\left(k+\frac{h}{q}+\beta\right)}{\left(\frac{h}{q}+\beta\right)_k}$$

$$= -\sum_{h=0}^{q-1} \frac{x^h}{\Gamma\left(\frac{h}{q}+\beta\right)}\left[\psi\left(\frac{h}{q}+\beta\right) + \sum_{k=1}^{\infty} \frac{(x^q)^k \psi\left(k+\frac{h}{q}+\beta\right)}{\left(\frac{h}{q}+\beta\right)_k}\right].$$
$$(5.6)$$

Finally, apply (2.36) to obtain (5.4), as we wanted to prove. □

**Definition 5.2** From (2.36), define the function:

$$P(a, x) := \frac{\partial Q(a, x)}{\partial x} = \psi(a)$$
$$+ e^x \left\{\frac{x-a-1}{a^2} {}_2F_2\left(\begin{array}{c}a, a\\a+1, a+1\end{array}\bigg| -x\right) + x^{-a}\gamma(a, x)\left[(x-a+1)\psi(a) + 1\right]\right\}.$$
$$(5.7)$$

**Theorem 5.3** *For $q \in \mathbb{N}$ and arbitrary $\beta$, the following reduction formula holds true:*

$$\left.\frac{\partial E_{\alpha,\beta}(x)}{\partial \alpha}\right|_{\alpha=1/q}$$
$$= -\sum_{h=0}^{q-1} \frac{h\left[\psi\left(\frac{h}{q}+\beta\right) + Q\left(\frac{h}{q}+\beta, x^q\right)\right] + q\,x^q P\left(\frac{h}{q}+\beta, x^q\right)}{\Gamma\left(\frac{h}{q}+\beta\right)} x^h. \quad (5.8)$$

**Proof** Apply (5.3) for $\alpha = 1/q$, taking into account (5.4),

$$\left.\frac{\partial E_{\alpha,\beta}(x)}{\partial \alpha}\right|_{\alpha=1/q} = x \frac{\partial}{\partial x}\left(\left.\frac{\partial E_{\alpha,\beta}(x)}{\partial \beta}\right|_{\alpha=1/q}\right)$$
$$= -\sum_{h=0}^{q-1} x \frac{\partial}{\partial x}\left(\frac{x^h}{\Gamma\left(\frac{h}{q}+\beta\right)}\left[\psi\left(\frac{h}{q}+\beta\right)+Q\left(\frac{h}{q}+\beta, x^q\right)\right]\right). \quad (5.9)$$

Apply (5.7) to obtain (5.8), as we wanted to prove. □

We present some examples of the reduction formulas given in (5.8) and (5.4) in Tables 6 and 7 respectively, using the help of MATHEMATICA program.

## 6 Differentiation of the Wright function with respect to parameters $\alpha$ and $\beta$

Direct differentiation of (1.2) with respect to $\alpha$ and $\beta$ gives

$$\frac{\partial W_{\alpha,\beta}(x)}{\partial \alpha} = -\sum_{k=1}^{\infty} \frac{\psi(\alpha k+\beta)}{k!\,\Gamma(\alpha k+\beta)} k\, x^k, \quad (6.1)$$

$$\frac{\partial W_{\alpha,\beta}(x)}{\partial \beta} = -\sum_{k=0}^{\infty} \frac{\psi(\alpha k+\beta)}{k!\,\Gamma(\alpha k+\beta)} x^k, \quad (6.2)$$

and, similar to (3.6), we have

$$\frac{\partial W_{\alpha,\beta}(x)}{\partial \alpha} = x \frac{\partial}{\partial x}\left(\frac{\partial W_{\alpha,\beta}(x)}{\partial \beta}\right). \quad (6.3)$$

We present next some reduction formulas of the derivative of the Wright function with respect to $\alpha$ and $\beta$ for $\alpha = 1$ and $\beta = 0, 1, 2$. Some of these results have been obtained in [3] with the aid of MATHEMATICA. Here, we derive analytically all the results. For this purpose, we use the sum formulas [4, Eqn. 6.2.1(19-20)]:

$$\sum_{k=0}^{\infty} \frac{x^k}{(k!)^2} \psi(k+1) = \frac{1}{2}\ln x\, I_0\left(2\sqrt{x}\right) + K_0\left(2\sqrt{x}\right), \quad (6.4)$$

and

$$\sum_{k=0}^{\infty} \frac{x^k}{k!\,(k+1)!} \psi(k+1) = \frac{2 - I_0\left(2\sqrt{x}\right) + \sqrt{x}\left[\ln x\, I_1\left(2\sqrt{x}\right) - 2K_1\left(2\sqrt{x}\right)\right]}{2x}.$$
$$(6.5)$$

**Table 6** First derivative of the Mittag-Leffler function $E_{\alpha,\beta}(x)$ with respect to parameter $\alpha$ derived for some values of parameters $\alpha$ and $\beta$

| $\alpha$ | $\beta$ | $\dfrac{\partial E_{\alpha,\beta}(x)}{\partial \alpha}$ |
|---|---|---|
| $\frac{1}{2}$ | $1$ | $-\dfrac{16x^5 e^{x^2} {}_2F_2\left(\frac{3}{2},\frac{3}{2};\frac{5}{2},\frac{5}{2};-x^2\right)}{9\sqrt{\pi}} - 2\left[e^{x^2}\left(1+2x^2\ln x - x^2\mathrm{Ei}\left(-x^2\right)\right)-1\right] - \dfrac{2}{\sqrt{\pi}}\left[\left(2x^3+x\right)\psi\left(\frac{3}{2}\right)+2e^{x^2}\gamma\left(\frac{3}{2},x^2\right)\left(1+x^2\psi\left(\frac{3}{2}\right)\right)\right]$ |
| $\frac{1}{2}$ | $2$ | $-\dfrac{4}{75\sqrt{\pi}x^2}\left[8\left(x^2-1\right)e^{x^2}x^5\,{}_2F_2\left(\frac{5}{2},\frac{5}{2};\frac{7}{2},\frac{7}{2};-x^2\right)+50e^{x^2}\gamma\left(\frac{5}{2},x^2\right)\left(x^2\psi\left(\frac{5}{2}\right)+1-\psi\left(\frac{5}{2}\right)\right)+25\left(2x^2+1\right)x^3\psi\left(\frac{5}{2}\right)\right]$ $-2x^2\left(\dfrac{e^{x^2}\left(\left(-\gamma x^2+x^2+\gamma\right)\gamma\left(2,x^2\right)+\left(1-x^2\right)\right)\left(-e^{-x^2}+\mathrm{Ei}(-x^2)-2\ln x-\gamma+1\right)}{x^4}-\gamma+1\right)$ |
| $\frac{1}{2}$ | $3$ | $-\dfrac{1}{x^4}\left[\left[\left(\frac{3}{2}-\gamma\right)\left(x^2-2\right)+1\right]e^{x^2}\gamma\left(3,x^2\right)-\left(x^2-2\right)\left[x^2+3+e^{x^2}\left(2\gamma+4\ln x-2\mathrm{Ei}\left(-x^2\right)-3\right)\right]\right]-x^2\left(\frac{3}{2}-\gamma\right)$ $-\dfrac{8}{735\sqrt{\pi}x^4}\left\{e^{x^2}\left[8\left(x^2-2\right)x^7\,{}_2F_2\left(\frac{7}{2},\frac{7}{2};\frac{9}{2},\frac{9}{2};-x^2\right)+98\gamma\left(\frac{7}{2},x^2\right)\left(\left(x^2-2\right)\psi\left(\frac{7}{2}\right)+1\right)\right]+49\left(2x^2+1\right)x^5\psi\left(\frac{7}{2}\right)\right\}$ |
| $\frac{1}{2}$ | $\beta$ | $-\dfrac{2x^2}{\Gamma(\beta)}\left\{e^{x^2}\left[\dfrac{\left(1-\beta+x^2\right){}_2F_2\left(\beta,\beta;\beta+1,\beta+1;-x^2\right)}{\beta^2}+x^{-2\beta}\gamma\left(\beta,x^2\right)\left[\psi(\beta)\left(1-\beta+x^2\right)+1\right]\right]+\psi(\beta)\right\}$ $-\dfrac{x}{\Gamma\left(\beta+\frac{1}{2}\right)}\left(\dfrac{8e^{x^2}x^2\left(1-\beta+x^2\right){}_2F_2\left(\beta+\frac{1}{2},\beta+\frac{1}{2};\beta+\frac{3}{2},\beta+\frac{3}{2};-x^2\right)}{(2\beta+1)^2}+\left(2x^2+1\right)\psi\left(\beta+\frac{1}{2}\right)+2e^{x^2}x^{1-2\beta}\gamma\left(\beta+\frac{1}{2},x^2\right)\left[\psi\left(\beta+\frac{1}{2}\right)\left(1-\beta+x^2\right)+1\right]\right)$ |
| $1$ | $\beta$ | $-\dfrac{x}{\Gamma(\beta)}\left\{e^x\left(\dfrac{(1-\beta+x){}_2F_2(\beta,\beta;\beta+1,\beta+1;-x)}{\beta^2}+x^{-\beta}\gamma(\beta,x)\left[\psi(\beta)(1-\beta+x)+1\right]\right)+\psi(\beta)\right\}$ |

**Table 7** First derivative of the Mittag-Leffler function $E_{\alpha,\beta}(x)$ with respect to parameter $\beta$ derived for some values of parameters $\alpha$ and $\beta$

| $\alpha$ | $\beta$ | $\dfrac{\partial E_{\alpha,\beta}(x)}{\partial \beta}$ |
|---|---|---|
| $\tfrac{1}{3}$ | 1 | $-\dfrac{9e^{x^3} x^4 \, _2F_2\left(\tfrac{4}{3},\tfrac{4}{3};\tfrac{7}{3},\tfrac{7}{3};-x^3\right)+16\psi\left(\tfrac{4}{3}\right)\left[x+e^{x^3}\gamma\left(\tfrac{4}{3},x^3\right)\right]}{16\,\Gamma\left(\tfrac{4}{3}\right)} - \dfrac{9e^{x^3} x^5 \, _2F_2\left(\tfrac{5}{3},\tfrac{5}{3};\tfrac{8}{3},\tfrac{8}{3};-x^3\right)+25\psi\left(\tfrac{5}{3}\right)\left[x^2+e^{x^3}\gamma\left(\tfrac{5}{3},x^3\right)\right]}{25\,\Gamma\left(\tfrac{4}{3}\right)} + e^{x^3}\left[\mathrm{Ei}\left(-x^3\right)-3\ln x\right]$ |
| $\tfrac{1}{3}$ | 2 | $-\dfrac{9e^{x^3} x^7 \, _2F_2\left(\tfrac{7}{3},\tfrac{7}{3};\tfrac{10}{3},\tfrac{10}{3};-x^3\right)+49\psi\left(\tfrac{7}{3}\right)\left[x^4+e^{x^3}\gamma\left(\tfrac{7}{3},x^3\right)\right]}{49\,x^3\,\Gamma\left(\tfrac{7}{3}\right)} - \dfrac{9e^{x^3} x^8 \, _2F_2\left(\tfrac{8}{3},\tfrac{8}{3};\tfrac{11}{3},\tfrac{11}{3};-x^3\right)+64\psi\left(\tfrac{8}{3}\right)\left[x^5+e^{x^3}\gamma\left(\tfrac{8}{3},x^3\right)\right]}{64\,x^3\,\Gamma\left(\tfrac{8}{3}\right)}$ $+ \dfrac{e^{x^3}\left\{\mathrm{Ei}\left(-x^3\right)-3\ln x+(\gamma-1)\left[\gamma\left(2,x^3\right)-1\right]\right\}-1}{x^3} + \gamma - 1$ |
| $\tfrac{1}{3}$ | 3 | $-\dfrac{9e^{x^3} x^{10} \, _2F_2\left(\tfrac{10}{3},\tfrac{10}{3};\tfrac{13}{3},\tfrac{13}{3};-x^3\right)+100\psi\left(\tfrac{10}{3}\right)\left[x^7+e^{x^3}\gamma\left(\tfrac{10}{3},x^3\right)\right]}{100\,x^6\,\Gamma\left(\tfrac{10}{3}\right)} - \dfrac{9e^{x^3} x^{11} \, _2F_2\left(\tfrac{11}{3},\tfrac{11}{3};\tfrac{14}{3},\tfrac{14}{3};-x^3\right)+121\psi\left(\tfrac{11}{3}\right)\left[x^8+e^{x^3}\gamma\left(\tfrac{11}{3},x^3\right)\right]}{121\,x^6\,\Gamma\left(\tfrac{11}{3}\right)}$ $+ \dfrac{e^{x^3}\left\{2\,\mathrm{Ei}\left(-x^3\right)-6\ln x+(2\gamma-3)\left[\tfrac{1}{2}\gamma\left(3,x^3\right)-1\right]\right\}-x^3-3}{2x^6} + \dfrac{\gamma}{2} - \dfrac{3}{4}$ |
| $\tfrac{1}{3}$ | $\beta$ | $-\dfrac{e^{x^3} x^3 \, _2F_2\left(\beta,\beta;\beta+1,\beta+1;-x^3\right)+\beta^2\psi(\beta)\left[1+e^{x^3}x^{3-3\beta}\gamma(\beta,x^3)\right]}{\beta\,\Gamma(\beta+1)} - \dfrac{x\left\{e^{x^3}x^3 \, _2F_2\left(\beta+\tfrac{1}{3},\beta+\tfrac{1}{3};\beta+\tfrac{4}{3},\beta+\tfrac{4}{3};-x^3\right)+\left(\beta+\tfrac{1}{3}\right)^2\psi\left(\beta+\tfrac{1}{3}\right)\left[1+e^{x^3}x^{2-3\beta}\gamma\left(\beta+\tfrac{1}{3},x^3\right)\right]\right\}}{\left(\beta+\tfrac{1}{3}\right)\Gamma\left(\beta+\tfrac{4}{3}\right)}$ |

**Table 7** continued

| $\alpha$ | $\beta$ | $\dfrac{\partial \mathrm{E}_{\alpha,\beta}(x)}{\partial \beta}$ |
|---|---|---|
| | | $-\dfrac{x^2\left\{e^{x^3}x^3\,_2F_2\left(\beta+\frac{2}{3},\beta+\frac{2}{3};\beta+\frac{5}{3},\beta+\frac{5}{3};-x^3\right)+\left(\beta+\frac{2}{3}\right)^2\psi\left(\beta+\frac{2}{3}\right)\left[1+e^{x^3}x^{1-3\beta}\gamma\left(\beta+\frac{2}{3},x^3\right)\right]\right\}}{\left(\beta+\frac{2}{3}\right)\Gamma\left(\beta+\frac{5}{3}\right)}$ |
| $\frac{1}{2}$ | 1 | $\frac{1}{9}e^{x^2}\left\{9\left[(\gamma-2+\ln 4)\mathrm{erf}\,x+\mathrm{Ei}\left(-x^2\right)-2\ln x\right]-\dfrac{8x^3}{\sqrt{\pi}}\,_2F_2\left(\tfrac{3}{2},\tfrac{3}{2};\tfrac{5}{2},\tfrac{5}{2};-x^2\right)\right\}$ |
| $\frac{1}{2}$ | 2 | $\dfrac{1}{75\sqrt{\pi}x^2}\left\{e^{x^2}\left[25\sqrt{\pi}\left((3\gamma-8+\ln 64)\mathrm{erf}\,x+3\mathrm{Ei}\left(-x^2\right)-6\ln x\right)-16x^5\,_2F_2\left(\tfrac{5}{2},\tfrac{5}{2};\tfrac{7}{2},\tfrac{7}{2};-x^2\right)\right]-25\left[3\gamma\left(2x+\sqrt{\pi}\right)+4x(\ln 8-4)\right]\right\}$ |
| $\frac{1}{2}$ | 3 | $\dfrac{1}{2205\sqrt{\pi}x^4}\left\{3e^{x^2}\left[49\sqrt{\pi}\left((15\gamma-46+15\ln 4)\mathrm{erf}\,x+15\left(\mathrm{Ei}\left(-x^2\right)-2\ln x\right)\right)-32x^7\,_2F_2\left(\tfrac{7}{2},\tfrac{7}{2};\tfrac{9}{2},\tfrac{9}{2};-x^2\right)\right]\right.$ $\left.-49\left[x\left(2\left(2x^2+3\right)(15\ln 4-46)-45\sqrt{\pi}x\right)+15\gamma\left(4x^3+3\sqrt{\pi}\left(x^2+1\right)+6x\right)\right]\right\}$ |

**Table 7** continued

| $\alpha$ | $\beta$ | $\dfrac{\partial E_{\alpha,\beta}(x)}{\partial \beta}$ |
|---|---|---|
| $\dfrac{1}{2}$ | $\beta$ | $-\dfrac{1}{\Gamma(\beta)}\left[e^{x^2}x^2\left(\dfrac{{}_2F_2(\beta,\beta;\beta+1,\beta+1;-x^2)}{\beta^2}+\psi(\beta)x^{-2\beta}\gamma(\beta,x^2)\right)+\psi(\beta)\right]$ |
|  | $\beta$ | $-\dfrac{x}{\Gamma\!\left(\beta+\tfrac{1}{2}\right)}\left[e^{x^2}x\left(\dfrac{x\,{}_2F_2\!\left(\beta+\tfrac{1}{2},\beta+\tfrac{1}{2};\beta+\tfrac{3}{2},\beta+\tfrac{3}{2};-x^2\right)}{\left(\beta+\tfrac{1}{2}\right)^2}+\psi\!\left(\beta+\tfrac{1}{2}\right)x^{-2\beta}\gamma\!\left(\beta+\tfrac{1}{2},x^2\right)\right)+\psi\!\left(\beta+\tfrac{1}{2}\right)\right]$ |
| $1$ | $\beta$ | $-\dfrac{1}{\Gamma(\beta)}\left[x\,e^x\left(\dfrac{{}_2F_2(\beta,\beta;\beta+1,\beta+1;-x)}{\beta^2}+\psi(\beta)\,x^{-\beta}\gamma(\beta,x)\right)+\psi(\beta)\right]$ |

**Theorem 6.1** *For $\alpha = 1$ and $\beta = 0$, the following reduction formula holds true:*

$$\left.\frac{\partial W_{\alpha,\beta}(x)}{\partial \beta}\right|_{\alpha=1, \beta=0} = \frac{1}{2}\left\{I_0\left(2\sqrt{x}\right) - \sqrt{x}\left[\ln x\, I_1\left(2\sqrt{x}\right) - 2K_1\left(2\sqrt{x}\right)\right]\right\} - 1. \tag{6.6}$$

*Proof* Note that

$$\left.\frac{\partial W_{\alpha,\beta}(x)}{\partial \beta}\right|_{\alpha=1, \beta=0} = -\sum_{k=0}^{\infty}\frac{x^k \psi(k)}{k!\,\Gamma(k)} = -\sum_{k=1}^{\infty}\frac{x^k \psi(k)}{k!\,(k-1)!} = -x\sum_{k=0}^{\infty}\frac{x^k \psi(k+1)}{(k+1)!\,k!}, \tag{6.7}$$

thus, applying (6.5), we obtain (6.6), as we wanted to prove. $\square$

**Theorem 6.2** *For $\alpha = 1$ and $\beta = 1$, the following reduction formula holds true:*

$$\left.\frac{\partial W_{\alpha,\beta}(x)}{\partial \beta}\right|_{\alpha=\beta=1} = -\frac{1}{2}\ln x\, I_0\left(2\sqrt{x}\right) - K_0\left(2\sqrt{x}\right). \tag{6.8}$$

*Proof* Directly from (6.4), we obtain the desired result

$$\left.\frac{\partial W_{\alpha,\beta}(x)}{\partial \beta}\right|_{\alpha=\beta=1} = -\sum_{k=0}^{\infty}\frac{x^k \psi(k+1)}{k!\,\Gamma(k+1)} = -\frac{1}{2}\ln x\, I_0\left(2\sqrt{x}\right) - K_0\left(2\sqrt{x}\right). \tag{6.9}$$

$\square$

**Theorem 6.3** *For $\alpha = 1$ and $\beta = 2$, the following reduction formula holds true:*

$$\left.\frac{\partial W_{\alpha,\beta}(x)}{\partial \beta}\right|_{\alpha=\beta=1} = \frac{2K_1\left(2\sqrt{x}\right) - \ln x\, I_1\left(2\sqrt{x}\right)}{2\sqrt{x}} - \frac{I_0\left(2\sqrt{x}\right)}{2x}. \tag{6.10}$$

*Proof* Apply (2.3) to obtain

$$\left.\frac{\partial W_{\alpha,\beta}(x)}{\partial \beta}\right|_{\alpha=1,\,\beta=2} = -\sum_{k=0}^{\infty}\frac{x^k \psi(k+2)}{k!\,\Gamma(k+2)} = -\sum_{k=0}^{\infty}\frac{x^k}{k!\,(k+1)!}\left[\frac{1}{k+1} + \psi(k+1)\right]$$

$$= -\sum_{k=0}^{\infty}\frac{x^k}{(k+1)!\,(k+1)!} - \sum_{k=0}^{\infty}\frac{x^k}{k!\,(k+1)!}\psi(k+1). \tag{6.11}$$

Note that, according to the definition of the modified Bessel function [16, Eqn. 10.25.2], we have

$$I_0\left(2\sqrt{x}\right) = \sum_{k=0}^{\infty} \frac{x^k}{k!k!} = 1 + \sum_{k=1}^{\infty} \frac{x^k}{k!k!} = 1 + x \sum_{k=0}^{\infty} \frac{x^k}{(k+1)!(k+1)!}. \quad (6.12)$$

Inserting (6.12) and (6.5) in (6.11) and grouping terms, we arrive at (6.10), as we wanted to prove. □

**Theorem 6.4** *For $\alpha = 1$ and $\beta = 0$, the following reduction formula holds true:*

$$\left.\frac{\partial W_{\alpha,\beta}(x)}{\partial \alpha}\right|_{\alpha=1,\,\beta=0} = -x\left[K_0\left(2\sqrt{x}\right) + \frac{1}{2}I_0\left(2\sqrt{x}\right)\ln x\right]. \quad (6.13)$$

*Proof* Apply (6.3) to the result given in (6.6). □

**Theorem 6.5** *For $\alpha = 1$ and $\beta = 1$, the following reduction formula holds true:*

$$\left.\frac{\partial W_{\alpha,\beta}(x)}{\partial \alpha}\right|_{\alpha=\beta=1} = \frac{\sqrt{x}\left[2K_1\left(2\sqrt{x}\right) - \ln x\, I_1\left(2\sqrt{x}\right)\right] - I_0\left(2\sqrt{x}\right)}{2}. \quad (6.14)$$

*Proof* Apply (6.3) to the result given in (6.8). □

**Theorem 6.6** *For $\alpha = 1$ and $\beta = 2$, the following reduction formula holds true:*

$$\left.\frac{\partial W_{\alpha,\beta}(x)}{\partial \alpha}\right|_{\alpha=1,\,\beta=2} = \frac{1}{2x}I_0\left(2\sqrt{x}\right) - \frac{1}{\sqrt{x}}I_1\left(2\sqrt{x}\right) - \frac{1}{2}\ln x\, I_2\left(2\sqrt{x}\right) - K_2\left(2\sqrt{x}\right). \quad (6.15)$$

*Proof* Apply (6.3) to the result given in (6.10). □

**Table 8** First derivative of the Wright function $W_{\alpha,\beta}(x)$ with respect to parameter $\alpha$ derived for some values of parameters $\alpha$ and $\beta$

| $\alpha$ | $\beta$ | $\frac{\partial W_{\alpha,\beta}(x)}{\partial \alpha}$ |
|---|---|---|
| 1 | 0 | $-x\left[K_0\left(2\sqrt{x}\right) + \frac{1}{2}I_0\left(2\sqrt{x}\right)\ln x\right]$ |
| 1 | 1 | $\frac{1}{2}\left\{\sqrt{x}\left[2K_1\left(2\sqrt{x}\right) - \ln x\, I_1\left(2\sqrt{x}\right)\right] - I_0\left(2\sqrt{x}\right)\right\}$ |
| 1 | 2 | $\frac{1}{2x}I_0\left(2\sqrt{x}\right) - \frac{1}{\sqrt{x}}I_1\left(2\sqrt{x}\right) - \frac{1}{2}\ln x\, I_2\left(2\sqrt{x}\right) - K_2\left(2\sqrt{x}\right)$ |

**Table 9** First derivative of the Wright function $W_{\alpha,\beta}(x)$ with respect to parameter $\beta$ derived for some values of parameters $\alpha$ and $\beta$

| $\alpha$ | $\beta$ | $\frac{\partial W_{\alpha,\beta}(x)}{\partial \beta}$ |
|---|---|---|
| 1 | 0 | $\frac{1}{2}\left\{I_0\left(2\sqrt{x}\right) - \sqrt{x}\left[\ln x\, I_1\left(2\sqrt{x}\right) - 2K_1\left(2\sqrt{x}\right)\right]\right\} - 1$ |
| 1 | 1 | $-\frac{1}{2}\ln x\, I_0\left(2\sqrt{x}\right) - K_0\left(2\sqrt{x}\right)$ |
| 1 | 2 | $\frac{2K_1(2\sqrt{x}) - \ln x\, I_1(2\sqrt{x})}{2\sqrt{x}} - \frac{I_0(2\sqrt{x})}{2x}$ |

We summarize the results of this section in Tables 8 and 9. Note that the series (6.1) and (6.2) are obtained in terms of the modified Bessel functions for $\alpha$ and $\beta$ being integers.

**Funding** Open Access funding provided thanks to the CRUE-CSIC agreement with Springer Nature.

## Declarations

**Conflict of interest** The authors declare that they have no conflict of interest

## References

1. Apelblat, A.: Differentiation of the Mittag-Leffler functions with respect to parameters in the Laplace transform approach. Mathematics **8**(5), 657 (2020)
2. Apelblat, A., González-Santander, J.L.: The Integral Mittag-Leffler. Whittaker and Wright functions. Mathematics **9**(24), 3255 (2021)
3. Apelblat, A., Mainardi, F.: Differentiation of the Wright functions with respect to parameters and other results. arXiv preprint arXiv:2009.08803 (2020)
4. Brychkov, Y.A.: Handbook of Special Functions: Derivatives, Integrals. Chapman and Hall/CRC, Series and Other Formulas (2008)
5. Chechkin, A., Gorenflo, R., Sokolov, I.: Retarding subdiffusion and accelerating superdiffusion governed by distributed-order fractional diffusion equations. Physical Review E **66**(4), 046129 (2002)
6. Chechkin, A.V., Klafter, J., Sokolov, I.M.: Fractional Fokker-Planck equation for ultraslow kinetics. Europhysics Letters **63**(3), 326 (2003)
7. Cvijović, D.: Closed-form summations of certain hypergeometric-type series containing the digamma function. Journal of Physics A: Mathematical and Theoretical **41**(45), 455205 (2008)
8. González-Santander, J.L., Sánchez Lasheras, F.: Finite and infinite hypergeometric sums involving the digamma function. Mathematics **10**(16), 2990 (2022)
9. Gorenflo, R., Kilbas, A.A., Mainardi, F., Rogosin, S.V., et al.: Mittag-Leffler Functions. Springer, Related Topics and Applications (2020)
10. Gorenflo, R., Luchko, Y., Mainardi, F.: Analytical properties and applications of the Wright functions. Fractional Calculus and Applied Analysis **2**(4), 383–414 (1999)

11. Górska, K., Pietrzak, T., Sandev, T., Tomovski, Z.: Volterra-Prabhakar derivative of distributed order and some applications. arXiv preprint arXiv:2212.13565 (2022)
12. Lebedev, N.: Special Functions and Their Applications. Prentice-Hall, Inc. (1965)
13. Mainardi, F.: Fractional Calculus and Waves in Linear Viscoelasticity: An Introduction to Mathematical Models, 2 edn. World Scientific (2022)
14. Miller, A.R.: Summations for certain series containing the digamma function. Journal of Physics A: Mathematical and General **39**(12), 3011 (2006)
15. Oldham, K.B., Myland, J., Spanier, J.: An Atlas of Functions: With Equator, the Atlas Function Calculator. Springer (2009)
16. Olver, F.W., Lozier, D.W., Boisvert, R.F., Clark, C.W.: NIST Handbook of Mathematical Functions. Cambridge University Press (2010)
17. Podlubny, I.: Fractional Differential Equations. Academic Press (1998)
18. Sandev, T., Tomovski, Z.: Langevin equation for a free particle driven by power law type of noises. Physics Letters A **378**(1–2), 1–9 (2014)
19. Sandev, T., Tomovski, Z., Dubbeldam, J.L., Chechkin, A.: Generalized diffusion-wave equation with memory kernel. Journal of Physics A: Mathematical and Theoretical **52**(1), 015201 (2018)
20. Wright, E.M.: On the coefficients of power series having exponential singularities. Journal of the London Mathematical Society **1**(1), 71–79 (1933)
21. Wright, E.M.: The generalized Bessel function of order greater than one. The Quarterly Journal of Mathematics **os–11**(1), 36–48 (1940)